\documentclass[final]{siamltex}
\usepackage{epsfig}
\usepackage{color}
\usepackage{amscd}
\usepackage{amsmath, amsfonts, amssymb,mathrsfs}
\usepackage{multirow}
\usepackage{enumerate}
\usepackage{verbatim}
\usepackage{cite}
\usepackage{tikz}
\usetikzlibrary{matrix,arrows}
\usepackage[justification=raggedright]{caption}
\usepackage{subcaption}
\usepackage[ruled]{algorithm2e}
\usepackage[export]{adjustbox}

\numberwithin{equation}{section}

\def \vec#1{{\bf{#1}}}
\def\pd#1#2{\frac{\partial #1}{\partial #2}}

\makeatletter
\def\hlinewd#1{%
\noalign{\ifnum0=`}\fi\hrule \@height #1 %
\futurelet\reserved@a\@xhline}
\makeatother

\newcommand{\bi}{\begin{itemize}}
\newcommand{\ei}{\end{itemize}}

\newcommand{\diverg}{\vec{\nabla}\cdot}
\newcommand{\director}{\vec{n}}
\newcommand{\curl}{\vec{\nabla}\times}
\newcommand{\Ltwoinner}[3]{\langle #1,#2 \rangle_0}
\newcommand{\Ltwonorm}[2]{\Vert #1 \Vert_0}
\newcommand{\Ltwonormndim}[3]{\Vert #1 \Vert_0}
\newcommand{\Ltwoinnerndim}[4]{\langle #1,#2 \rangle_0}
\newcommand{\Rthree}{\mathbb{R}^3}

\newcommand{\diff}[1]{\, d#1}

\newcommand{\kdirector}{\vec{n}_k}
\newcommand{\ddirector}{\delta \director}
\newcommand{\dlambda}{\delta \lambda}
\newcommand{\klambda}{\lambda_k}

\newcommand{\kphi}{\phi_k}

\newcommand{\dphi}{\delta \phi}
\newcommand{\lagdivn}{\mathcal{L}_{\director}[\vec{v}]}
\newcommand{\lagdivlam}{\mathcal{L}_{\lambda}[\gamma]}

\newcommand{\lagdivphi}{\mathcal{L}_{\phi}[\psi]}

\newcommand{\Honenorm}[2]{\Vert #1 \Vert_1}

\newcommand{\Hdc}{\mathcal{H}^{DC}{(\Omega)}}
\newcommand{\Hdcnot}{\mathcal{H}^{DC}_0{(\Omega)}}

\newcommand{\Hone}[1]{H^1(#1)}

\newcommand{\Honenot}[1]{H^1_0({#1})}
\newcommand{\Honebnot}[1]{H^{1,0}({#1})}
\newcommand{\Honeb}[1]{H^{1,g}(#1)}

\newcommand{\Ltwo}[1]{L^2(#1)}

\newcommand{\Lp}[1]{L^p (\Omega)}
\newcommand{\Unorm}[1]{\Vert #1 \Vert_{U}}
\newcommand{\Hcurl}[1]{H(\text{curl},#1)}
\newcommand{\Hdiv}[1]{H(\text{div},#1)}

\newcommand{\Linfinity}[1]{L^{\infty}(\Omega)}

\newcommand{\Hdcnorm}[2]{\Vert #1 \Vert_{DC}}

\newcounter{casenum}

\title{A Deflation Technique for Detecting Multiple Liquid Crystal Equilibrium States\thanks{Submitted \today}}

\author{J. H. Adler\thanks{Department of Mathematics, Tufts University, Medford, MA 02155 (james.adler@tufts.edu, david.emerson@tufts.edu).} \and D. B. Emerson\footnotemark[2] \and P. E. Farrell\thanks{Mathematical Institute, University of Oxford, Oxford OX2 6GG, United Kingdom (patrick.farrell@maths.ox.ac.uk). This author's work was funded by EPSRC grants EP/K030930/1, EP/M011151/1, and a Center of Excellence grant from the Research Council of Norway to the Center for
Biomedical Computing at Simula Research Laboratory.} \and S. P. MacLachlan\thanks{Department of Mathematics and Statistics, Memorial University of Newfoundland, St. John's, Newfoundland and Labrador A1C 5S7, Canada (smaclachlan@mun.ca). This author's work was partially supported by an NSERC Discovery Grant.}}

\begin{document}

\maketitle

\begin{abstract}
Multiple equilibrium states arise in many physical systems, including
various types of liquid crystal structures. Having the ability to
reliably compute such states enables more accurate physical analysis
and understanding of experimental behavior. This paper adapts and
extends a deflation technique for the computation of multiple distinct
solutions arising in the context of modeling equilibrium
configurations of nematic and cholesteric liquid crystals. The
deflation method is applied as part of an overall free-energy
variational approach and is modified to fit the framework of
optimization of a functional with pointwise constraints. It is shown that multigrid methods designed for the undeflated systems may be applied to efficiently solve the linear systems arising in the application of deflation. For the numerical algorithm, the deflation approach is interwoven with nested iteration, creating a dynamic and efficient method that further enables the discovery of distinct solutions. Finally, four numerical experiments are performed demonstrating the efficacy and accuracy of the algorithm in detecting important physical phenomena, including bifurcation and disclination behaviors. The final numerical experiment expands the algorithm to model cholesteric liquid crystals and illustrates the full discovery power of the deflation process.
\end{abstract}


 
\pagestyle{myheadings}
\thispagestyle{plain}
\markboth{\sc Adler, Emerson, Farrell, Maclachlan}{\sc Deflation for Liquid Crystal Equilibria}


\section{Introduction}

As materials with mesophases exhibiting characteristics of both liquids and structured solids, liquid crystals produce a striking variety of arrangements and behaviors. These mesophases are found at varying temperatures and solvent concentrations and exist for many types of materials, including both synthetic \cite{Gattermann1} and naturally occurring molecular compositions \cite{Collings2}. In this paper, we focus on nematic phases, which consist of rod-like molecules, and cholesteric liquid crystals, which share many similarities with nematics but intrinsically prefer helical structures that admit less symmetry due to chiral preference. These types of liquid crystals self-assemble into ordered structures characterized by a preferred average direction at each point known as the director. The director is described by a unit vector field at each point and is denoted $\director(x, y, z) = (n_1(x, y, z), n_2(x, y, z), n_3(x, y, z))^T$. 

Along with their crystalline self-structuring, liquid crystals demonstrate a number of important physical phenomena including birefringence, electric coupling, and flexoelectric effects. Comprehensive reviews of liquid crystal physics are found in \cite{Stewart1, deGennes1, Virga1}. These properties and others have led to many important discoveries and a diversity of applications. In addition to their value for display technologies, liquid crystals are used in photorefractive cells \cite{Lagerwall1}, chemical sensing \cite{Shah1}, and the production of novel actuators \cite{Yamada1}. Numerical simulations of liquid crystal configurations are used to test and examine theory, explore new physical phenomena \cite{Atherton1, Emerson2}, and optimize device performance. Many modern experimental designs and physical effects, including bistability \cite{Stalder1}, require accurate and efficient large-scale simulations.

Herein, we consider the Frank-Oseen free-energy model for the computation of liquid crystal equilibrium configurations \cite{Stewart1, Virga1}. The complexity of the model and the necessary nonlinear pointwise unit-length constraint have limited the availability of analytical solutions in the absence of significant simplifying assumptions. Recently, a number of numerical methods \cite{Ramage1, Ramage3, Atherton1, Pandolfi1} have been developed for the Frank-Oseen model. In \cite{Emerson1, Emerson2}, a theoretically supported energy-minimization finite-element technique was developed that facilitates highly accurate and efficient computational simulation of complicated liquid crystal behavior. However, the presence of multiple local extrema and saddle-points can increase the difficulty of locating global extrema with existing methods.

In the context of static liquid crystal structures, it is well known that the associated system of partial differential equations (PDEs), commonly referred to as the Equilibrium Equations \cite{Ericksen4, Stewart1}, permit multiple solutions, even under relatively mild complexity \cite{Deuling1}. Further, multiple local minimizers and saddle-point solutions may exist in the energy formulation. In addition to increasing confidence in, and facilitating the computation of, global extrema, locating these distinct solutions reveals configurations with physical relevance, as is the case, for instance, with defect arrangements. In this paper, we adapt and expand the deflation methodology first proposed in \cite{Farrell1}. It has been successfully applied to compute distinct solutions to nonlinear PDEs and certain kinds of complementarity problems \cite{Farrell1, Farrell3}. 

The applied deflation technique sequentially modifies the nonlinear problem
to eliminate previously known solutions from
consideration, allowing for successive discovery of distinct solutions
to the nonlinear system under consideration. Here, we examine the
method's performance in the context of functional optimization subject
to a pointwise (nonlinear) constraint. The deflation method proves
particularly attractive as it allows for the preservation of
finite-element sparsity, admits the use of existing, advanced
multigrid methods, and seamlessly integrates with efficient nested
iteration routines. Moreover, the theory developed in \cite{Farrell1}
surrounding the deflation operators is directly applicable to the
spaces and operators used in this manuscript. Finally, the approach
is, in practice, highly successful in locating multiple, distinct equilibrium states.

This paper is organized as follows. The energy model and undeflated minimization approach are summarized in Section \ref{Model}. The deflation technique is discussed and derived in the context of the minimization framework in Section \ref{DeflationMethod}. The integration of previously designed multigrid methods and the interweave of deflation with nested iterations is also examined in the section. In Sections \ref{NumericalExperiments} and \ref{CholestericNumericalExperiments}, the algorithm implementation is outlined and four numerical experiments are performed. Finally, Section \ref{Conclusions} provides some concluding remarks and a discussion of future work.

\section{Nematic Energy Model and Minimization} \label{Model}

While a number of liquid crystal models exist \cite{Onsager1, Davis1,
  Stewart1}, we consider the Frank-Oseen free-energy model where the
equilibrium free energy for a domain $\Omega$ is characterized by deformations of the
nondimensional unit-length director field, $\director$. Liquid crystals tend towards configurations exhibiting minimal free
energy. Let $K_i$, $i = 1, 2, 3$, be the Frank constants \cite{Frank1}
with $K_i \geq 0$ \cite{Ericksen2}. Herein, we consider the case that each
$K_i \neq 0$.  These constants are often anisotropic (i.e., $K_1 \neq K_2 \neq K_3$), vary with liquid crystal type, and play important roles in physical phenomena \cite{Atherton2, Lee1}. 

We denote the classical $\Ltwo{\Omega}$ inner product and norm as $\Ltwoinner{\cdot}{\cdot}{\Omega}$ and $\Ltwonorm{\cdot}{\Omega}$, respectively, for both scalar and vector quantities. Throughout this paper, we assume the presence of Dirichlet boundary conditions or mixed Dirichlet and periodic boundary conditions on a rectangular domain and, therefore, utilize the null Lagrangian simplification discussed in \cite{Emerson2, Stewart1}. Hence, including the possibility of external electric fields, the Frank-Oseen free energy for nematics is written
\begin{align} \label{FrankOseenFree}
\int_{\Omega} \left( w_F - \frac{1}{2} \vec{D} \cdot \vec{E} \right ) \diff{V} &= \frac{1}{2} K_1 \Ltwonorm{\nabla\cdot \director}{\Omega}^2 + \frac{1}{2} K_3\Ltwoinnerndim{\vec{Z} \nabla \times \director}{\nabla \times \director}{\Omega}{3} \nonumber \\
& \quad - \frac{1}{2} \epsilon_0 \epsilon_{\perp}\Ltwoinnerndim{\nabla \phi}{\nabla\phi}{\Omega}{3}  - \frac{1}{2} \epsilon_0 \epsilon_a \Ltwoinner{\director \cdot \nabla\phi}{\director \cdot \nabla\phi}{\Omega},
\end{align}
where $\phi$ is an electric potential such that $\vec{E} =
-\nabla\phi$, $\epsilon_0$ denotes the
permittivity of free space, and the dimensionless constants
$\epsilon_{\perp}$ and $\epsilon_a$ are the perpendicular dielectric
permittivity and dielectric anisotropy of the liquid crystal,
respectively.  Finally, $\vec{Z} = \kappa \director \otimes \director
+ (\vec{I} - \director \otimes \director) = \vec{I} - (1- \kappa)
\director \otimes \director$, is a dimensionless tensor, 
where $\kappa = K_2/K_3$. Note that if $\kappa =1$, $\vec{Z}$ is reduced to the identity. 

In order to properly formulate the Lagrangian, a nondimensionalization
was introduced in \cite{Emerson3}, using a characteristic length
scale, $\sigma$, characteristic Frank constant, $K$, and
characteristic voltage, $\phi_0 > 0$. We apply a spatial change of
variables to the Frank-Oseen free energy along with a simplification, so that the entire
expression is dimensionless, along with the Frank constants, electric
potential, and the parameter $\epsilon_0$.  Note that this scaling
implies that the domain, $\Omega$, and differential operator,
$\nabla$, are also dimensionless.

As noted above, the director field is subject to a local unit-length constraint such that $\director \cdot \director = 1$ at each point throughout the domain. In \cite{Emerson3}, numerical evidence suggests that imposing this constraint with Lagrange multipliers is an accurate and highly efficient approach, particularly in comparison to penalty or renormalization formulations. Thus, to compute free-energy minimizing configurations, we define the nondimensionalized nematic free-energy functional, after rescaling by a factor of $2$, as
\begin{align}
\mathcal{F}(\director, \phi) &= K_1 \Ltwonorm{\diverg \director}{\Omega}^2 + K_3\Ltwoinnerndim{\vec{Z} \curl \director}{\curl \director}{\Omega}{3} - \epsilon_0 \epsilon_{\perp} \Ltwoinnerndim{\nabla \phi}{\nabla \phi}{\Omega}{3} \nonumber \\
&\qquad - \epsilon_0 \epsilon_a \Ltwoinner{\director \cdot \nabla \phi}{\director \cdot \nabla \phi}{\Omega}. \label{NematicFunctional}
\end{align}

Throughout this paper, we will make use of the spaces $\Hdiv{\Omega} = \{\vec{v} \in \left (L^2(\Omega) \right)^3 : \diverg \vec{v} \in L^2(\Omega) \}$ and 
$\Hcurl{\Omega} = \{ \vec{v} \in \left (L^2(\Omega) \right)^3 : \curl \vec{v} \in \left(L^2(\Omega) \right)^3 \}$. As in \cite{Emerson2}, define 
\begin{equation*}
\Hdc= \{ \vec{v} \in \Hdiv{\Omega} \cap \Hcurl{\Omega} : B(\vec{v}) = \bar{\vec{g}} \},
\end{equation*}
with norm $\Hdcnorm{\vec{v}}{\Omega}^2 =
\Ltwonormndim{\vec{v}}{\Omega}{3}^2 + \Ltwonorm{\diverg
  \vec{v}}{\Omega}^2 + \Ltwonormndim{\curl \vec{v}}{\Omega}{3}^2$ and
appropriate boundary conditions $B(\vec{v})=\bar{\vec{g}}$. Here, we
assume that $\bar{\vec{g}}$ satisfies appropriate compatibility
conditions for the operator $B$. For example, if $B$ represents full
Dirichlet boundary conditions and $\Omega$ has a Lipschitz continuous
boundary, it is assumed that $\bar{\vec{g}} \in
H^{\frac{1}{2}}(\partial \Omega)^3$ \cite{Girault1}. Further, let
$\Hdcnot = \{ \vec{v} \in \Hdiv{\Omega} \cap \Hcurl{\Omega} :
B(\vec{v}) = \vec{0} \}$. Note that if $\Omega$ is a Lipshitz domain
and $B$ imposes full Dirichlet boundary conditions on all components
of $\vec{v}$, then $\Hdcnot = \left (\Honenot{\Omega}\right )^3$ \cite[Lemma 2.5]{Girault1}. Denote
\begin{align*}
\Honeb{\Omega} = \{f \in \Hone{\Omega} : B_1(f) = g \},
\end{align*}
where $\Hone{\Omega}$ represents the classical Sobolev space with norm $\Honenorm{\cdot}{\Omega}$ and $B_1(f) = g$ is an appropriate boundary condition expression for $\phi$.

We define the Lagrangian as
\begin{align*}
\mathcal{L}(\director, \phi, \lambda) &= \mathcal{F}(\director, \phi) + \int_{\Omega} \lambda(\vec{x})(\director \cdot \director-1) \diff{V},
\end{align*}
where $\mathcal{L}(\director, \phi, \lambda)$ has been nondimensionalized in the same fashion as the free-energy functional. To minimize the functional, first-order optimality conditions are derived as
\begin{align}
\lagdivn &= \frac{\partial}{\partial \director} \mathcal{L}(\director, \phi, \lambda) [\vec{v}] =0 & & \forall \vec{v} \in \Hdcnot,  \label{FOOC1} \\
\lagdivphi &= \frac{\partial}{\partial \phi} \mathcal{L}(\director, \phi, \lambda) [\psi] =0 & & \forall \psi \in \Honebnot{\Omega}, \\
\lagdivlam &= \frac{\partial}{\partial \lambda} \mathcal{L}(\director, \phi, \lambda) [\gamma] =0 & & \forall \gamma \in L^2(\Omega) \label{FOOC3}.
\end{align}
Define the product space $U = \Hdc \times \Honeb{\Omega} \times \Ltwo{\Omega}$ with associated norm $\Unorm{\cdot}$, and denote the subspace $U_0 = \Hdcnot \times \Honebnot{\Omega} \times \Ltwo{\Omega}$. Further, let $\mathcal{A}(\director, \phi, \lambda; \vec{v}, \psi, \gamma): U \to \Rthree$ be the variational system operator for variations $\vec{v}$, $\psi$, and $\gamma$. The operator is expressed in component form as 
\begin{align*}
\mathcal{A}(\director, \phi, \lambda; \vec{v}, \psi, \gamma) = \left [\text{ } \lagdivn \text{ } \lagdivphi \text{ } \lagdivlam \text{ } \right]^T. 
\end{align*}
Thus, Equations \eqref{FOOC1}-\eqref{FOOC3} are more compactly written as seeking $(\director, \phi, \lambda) \in U$ such that
\begin{align}
\mathcal{A}(\director, \phi, \lambda; \vec{v}, \psi, \gamma) = \vec{0} & & \forall (\vec{v}, \psi, \gamma) \in U_0. \label{compactFOOC}
\end{align}
The above variational system is nonlinear and, under certain conditions, admits several distinct solutions. For example, the classical Freedericksz transition problem \cite{Zocher1, Freedericksz1}, which is discussed in detail below, admits at least three solutions to the first-order optimality conditions. 

\section{Deflation Methodology} \label{DeflationMethod}

In \cite{Emerson1, Emerson2}, Newton linearizations and finite
elements are used to compute solutions to the variational system
described in \eqref{compactFOOC}. A standard but unsystematic approach
to computing distinct solutions for nonlinear problems with several
solutions is the use of numerous initial guesses as part of an
overarching Newton-type scheme, known as multistart methods
\cite{Marti1}. In this section, we adapt the deflation technique first
proposed in \cite{Farrell1} as a more effective and systematic
alternative. Due to the existence of multiple solutions to the
variational form of the Euler-Lagrange equations in
\eqref{compactFOOC}, the question of whether presented solutions
represent global minima or only local minima (or maxima or saddle
points) is often difficult to answer with certainty. Furthermore,
under certain conditions, local extrema are observable and represent
valid physical states of a device in equilibrium. The deflation technique presented in this section systematically promotes the discovery of numerous equilibrium points, increasing confidence in global minimizer claims and revealing physically stable local minima.

Linearizing the undeflated variational system in \eqref{compactFOOC} yields the Newton update equations
\begin{equation} \label{ElectricPotentialNewtonSystem}
\left [ \begin{array}{c c c}
\mathcal{L}_{\director \director} & \mathcal{L}_{\director \phi} & \mathcal{L}_{\director \lambda} \\
\mathcal{L}_{\phi \director} & \mathcal{L}_{\phi \phi} & \vec{0} \\
\mathcal{L}_{\lambda \director} & \vec{0} & \vec{0}
\end{array} \right] 
\left [ \begin{array}{c}
\ddirector \\
\dphi \\
\dlambda
\end{array} \right] = -
\left[ \begin{array}{c}
\mathcal{L}_{\director} \\
\mathcal{L}_{\phi} \\
\mathcal{L}_{\lambda} 
\end{array} \right],
\end{equation}
where each of the system components is evaluated at $\kdirector$, $\kphi$, and $\lambda_k$, the current approximations for $\director$, $\phi$, and $\lambda$, and $\ddirector= \director_{k+1} - \kdirector$, $\dphi = \phi_{k+1} - \kphi$, and $\dlambda = \lambda_{k+1}-\klambda$ are the updates we seek to compute. The matrix-vector multiplication indicates the direction that the derivatives in the Hessian are taken. For instance, $\mathcal{L}_{\lambda \director}[\gamma] \cdot \ddirector = \pd{ }{\director} \left( \mathcal{L}_{\lambda} (\kdirector, \kphi, \klambda)[\gamma] \right)[\ddirector]$, where the partials indicate G\^{a}teaux derivatives in the respective variables. The complete system is found in \cite{Emerson2}.

Denote the identity matrix of appropriate dimension as $\vec{I}$, and let $\vec{r} = (\director_r, \phi_r, \lambda_r)$ represent a known solution to \eqref{compactFOOC}. Further, let $\vec{u} = (\director, \phi, \lambda)$. We define a shifted deflation operator,
\begin{align*}
M_{p, \alpha}(\vec{u}; \vec{r}) = \left (\frac{1}{\Unorm{\vec{u}-\vec{r}}^p} + \alpha \right) \vec{I},
\end{align*}
where $\alpha \geq 0$ is a shift scalar and $p \in [1, \infty)$ is the deflation exponent. Note that for a given $\vec{u}$ and $\vec{r}$, the deflation operator $M_{p, \alpha}(\vec{u}; \vec{r}): \Rthree \to \Rthree$. Applying the deflation operator to the variational operator $\mathcal{A}$ reduces the solution space by eliminating $\vec{r}$ as a possible zero of the deflated system. The resulting deflated variational operator is given by
\begin{align*}
\mathcal{G}(\director, \phi, \lambda; \vec{v}, \psi, \gamma) = M_{p, \alpha}(\vec{u}; \vec{r}) \mathcal{A}(\director, \phi, \lambda; \vec{v}, \psi, \gamma)
&= \left (\frac{1}{\Unorm{\vec{u}-\vec{r}}^p} + \alpha \right) \left [\begin{array}{c} \lagdivn \\ \lagdivphi \\ \lagdivlam \end{array} \right].
\end{align*}
This produces the deflated variational system
\begin{align}
\mathcal{G}(\director, \phi, \lambda; \vec{v}, \psi, \gamma) = \vec{0} & & \forall (\vec{v}, \psi, \gamma) \in U_0. \label{deflatedVarSystem}
\end{align}
The shift is used so that the deflated residual does not tend to zero as $\Unorm{\vec{u}-\vec{r}}$ becomes arbitrarily large, see \cite{Farrell1}. While the method is generally robust with respect to parameter choice, there are situations where additional performance improvements are attainable for certain selections of $p$ and $\alpha$. For brevity, we suppress the semicolon notation in the variational operators except when necessary for clarity and denote $\eta(\vec{u}) = \left (\frac{1}{\Unorm{\vec{u}-\vec{r}}^p} + \alpha \right)$. Note that the deflated variational operator, $\mathcal{G}(\vec{u}) = \eta(\vec{u}) \mathcal{A}(\vec{u})$, is also nonlinear. As a vector-valued functional, the linearization of $\mathcal{G}(\vec{u})$ requires computation of the Jacobian,
\begin{align}
J(\mathcal{G}(\vec{u}_k))[\delta \vec{u}] = \eta(\vec{u}_k) J(\mathcal{A}(\vec{u}_k))[\delta \vec{u}] + \mathcal{A}(\vec{u}_k) \otimes \nabla \eta(\vec{u}_k)[\delta \vec{u}], \label{DeflatedJacobianDef}
\end{align}
where $J(\mathcal{A}(\vec{u}_k))[\delta \vec{u}]$ represents the
Jacobian of $\mathcal{A}$ in the directions $\delta \vec{u} =
(\ddirector, \dphi, \dlambda)$,  $\nabla \eta(\vec{u}_k)[\delta
\vec{u}]$ denotes the gradient of $\eta$, with each evaluated at
$\vec{u}_k$, and $\otimes$ denotes the standard outer product of two vectors. Computing the Jacobian of $\mathcal{A}$ is equivalent to deriving the Hessian in \eqref{ElectricPotentialNewtonSystem}, previously computed in \cite{Emerson2}. Thus, all that is left to compute is the gradient of $\eta$. This gradient has the form
\begin{align} \label{EtaGradient}
\nabla \eta(\vec{u}_k)[\delta \vec{u}] = -\frac{p}{2}\Unorm{\vec{u}_k - r}^{-p-2} \left [ 
\begin{array}{c} 
\pd{}{\director} (\Unorm{\vec{u}_k-r}^2) [\ddirector] \\
 \pd{}{\phi} (\Unorm{\vec{u}_k-r}^2) [\dphi] \\
  \pd{}{\lambda} (\Unorm{\vec{u}_k-r}^2) [\dlambda]
 \end{array} \right ].
\end{align}
The G\^{a}teaux derivatives in \eqref{EtaGradient} expand to
\begin{align*}
\left [ 
\begin{array}{c} 
2 \Ltwoinnerndim{\ddirector}{\kdirector - \director_r}{\Omega}{3} + 2 \Ltwoinner{\diverg \ddirector}{\diverg (\kdirector - \director_r)}{\Omega} + 2 \Ltwoinnerndim{\curl \ddirector}{\curl (\kdirector - \director_r)}{\Omega}{3} \\
2 \Ltwoinner{\dphi}{\kphi - \phi_r}{\Omega} + 2 \Ltwoinnerndim{\nabla \dphi}{\nabla (\kphi - \phi_r)}{\Omega}{3}\\
2 \Ltwoinner{\dlambda}{\klambda - \lambda_r}{\Omega}
 \end{array} \right ].
\end{align*}
Constructing the Jacobian with the gradient in \eqref{EtaGradient}, the linearized system for the deflated problem is summarized as
\begin{align}
J(\mathcal{G}(\vec{u}_k))[\delta \vec{u}] = - \mathcal{G}(\vec{u}_k) & & \forall (\vec{v}, \psi, \gamma) \in U_0. \label{DeflatedSystem}
\end{align}

\subsection{Deflated Linear Systems} \label{DeflationMatrixFree}

In this section, we consider the structure of the linear systems arising from a mixed finite-element discretization \cite{Emerson2} of the linearized deflation problem in \eqref{DeflatedSystem}. Let $A(\vec{u}_k)$ and $G(\vec{u}_k)$ denote the vectors corresponding to discretizations of $\mathcal{A}(\vec{u}_k)$ and $\mathcal{G}(\vec{u}_k)$, respectively, and let $d(\vec{u}_k)$ be the discretization vector corresponding to the gradient of $\eta$. Let $J_G(\vec{u}_k)$ and $J_A(\vec{u}_k)$ indicate the discretized Jacobians of the deflated and undeflated systems, respectively. Then,
\begin{align*}
J_G(\vec{u}_k) = \eta(\vec{u}_k) J_A(\vec{u}_k) + A(\vec{u}_k) d(\vec{u}_k)^T.
\end{align*}
As defined, $J_G(\vec{u}_k)$ is composed of a rank-one update to $J_A(\vec{u}_k)$. Thus, $J_G(\vec{u}_k)$ is generally dense even if $J_A(\vec{u}_k)$ is not and explicit construction and computation with the matrix is prohibitively expensive. 

A strategy for constructing effective preconditioners for the deflated system based on existing preconditioners for the undeflated matrices and computing their actions in a matrix-free fashion is presented in \cite{Farrell1}. 
Here, we are interested in the reuse of fast iterative solvers designed for the original linear system. Specifically, we consider applying the coupled multigrid method with Braess-Sarazin-type relaxation developed in \cite{Emerson3} to solve for the Newton updates in the deflated linear system. 

Denote the discretization of the system right-hand-side in \eqref{DeflatedSystem} as $b_G(\vec{u}_k) = -\eta(\vec{u}_k) A(\vec{u}_k)$. Throughout the remainder of the paper, except when necessary for clarity, we neglect the dependence on $\vec{u}_k$ in the notation. This yields the compact representation $J_G = (\eta J_A + Ad^T)$. Applying the Sherman-Morrison formula \cite{Hager1}
\begin{align} \label{InverseJG}
J_G^{-1} = (\eta J_A + Ad^T)^{-1} = \frac{J_A^{-1}}{\eta} - \frac{\frac{1}{\eta^2}J_A^{-1}Ad^T J_A^{-1}}{1+ \frac{1}{\eta}d^T J_A^{-1}A}.
\end{align}
Using \eqref{InverseJG} to compute the update vector produces
\begin{align*}
J_G^{-1} b_G = \frac{J_A^{-1}b_G}{\eta} - \frac{\frac{1}{\eta^2}J_A^{-1} A d^T J_A^{-1} b_G}{1 + \frac{1}{\eta} d^T J_A^{-1} A} &= -J_A^{-1}A - \frac{-\frac{1}{\eta}J_A^{-1} A d^T J_A^{-1} A}{1 + \frac{1}{\eta} \cdot d^T J_A^{-1}A} \\
&= -\left(1 - \frac{\frac{1}{\eta} \cdot d^T J_A^{-1} A}{1 + \frac{1}{\eta} \cdot d^T J_A^{-1}A} \right) J_A^{-1}A.
\end{align*}
Note that $J_A^{-1}A$ corresponds to assembling and solving the
original undeflated problem and $d^TJ_A^{-1}A$ is a dot product
resulting in a scalar. Thus, solving the discrete form of the
deflation system in \eqref{DeflatedSystem} is reduced to a single
solve with the original sparse system, one dot product, one vector
scaling, and a few scalar operations. Therefore, the coupled multigrid
method from \cite{Emerson3} can be directly applied to the discrete system, $J_A$, to efficiently compute both deflated and undeflated Newton updates. 

\subsection{Multiple Deflation}

Thus far, the class of deflation operators considered focuses on deflation with one known solution, $\vec{r} = (\director_r, \phi_r, \lambda_r)$. In this section, we briefly discuss extending the deflation procedure to treat a family of known solutions $\vec{r}_1, \vec{r}_2, \ldots, \vec{r}_m$. With several known solutions, the multiple deflation operator is the product of the single deflation operators for each individual solution such that
\begin{align*}
M_{p, \alpha}(\vec{u}; \vec{r}_1, \vec{r}_2, \ldots, \vec{r}_m) = \prod_{i=1}^m M_{p, \alpha} (\vec{u}; \vec{r}_i).
\end{align*}
This modifies the action of $M_{p, \alpha}(\vec{u}, \vec{r}_1, \vec{r}_2, \ldots, \vec{r}_m)$ on $\mathcal{A}$ such that
\begin{align*}
\mathcal{G}(\director, \phi, \lambda; \vec{v}, \psi, \gamma) = \prod_{i=1}^m M_{p, \alpha}(\vec{u}, \vec{r}_i) \mathcal{A}(\director, \phi, \lambda) &= \left(\prod_{i=1}^m \left( \frac{1}{\Unorm{\vec{u}-\vec{r}_i}^p} + \alpha \right) \right) \left [\begin{array}{c} \lagdivn \\ \lagdivphi \\ \lagdivlam \end{array} \right], \\
\end{align*}
which we recognize, as in the case of single deflation, to be of the
form $\mathcal{G}(\vec{u}) = \eta(\vec{u})\mathcal{A}(\vec{u})$.
This deflated system remains nonlinear and corresponding
linearizations are derived to compute distinct solutions satisfying
the first-order optimality conditions. As with the single deflation
linearization, the multiple deflation Jacobian, $J_G$, is
composed of a rank-one update to $J_A$ as in
\eqref{DeflatedJacobianDef}, though $\nabla \eta$ is now more
complicated than the single deflation case in \eqref{EtaGradient}. A process similar to that applied in the single deflation case reveals an analogous result for computation of solutions to the discretized, deflated linearizations and yields similar results enabling the application of multigrid methods to linear systems subject to deflation over several known solutions. Each of the simulations to follow employs multi-solution deflation operators as distinct solutions are discovered.

\subsection{Interaction with Nested Iteration}

Nested iteration (NI) is a common tool for the numerical solution of
nonlinear PDEs \cite{Starke1}, where the system is first solved on a very coarse
level, where computation is cheap.  A series of refinement steps are
then taken, interpolating the coarse-grid solution to a finer mesh and
using this as an initial guess for the fine-grid problem.  A key
advantage is that these interpolated approximations are typically very
good initial guesses for Newton's method on the finer grids, so very
few iterations are needed on fine levels, where computation is expensive.
Such an NI process readily admits integration with the deflation
methodology through a combination of
continued iteration on known solutions on each level, followed by applying deflation
to uncover additional solutions on each mesh in the grid ordering. The
general numerical flow is detailed below in Algorithm $1$. The
algorithm has four main stages. The outermost phase is a nested
iteration hierarchy that has proven highly effective in
reducing computational work for these types of problems \cite{Emerson1,
  Emerson3}. On each mesh, the algorithm first performs (undeflated)
Newton iterations on interpolated versions of the solutions found on
the previous, coarser mesh, termed the \emph{continuation list} in Algorithm
\ref{algo}, to further resolve the solution features on the finer
mesh. This procedure is followed by a solution discovery stage
incorporating the set of known solutions to form deflated systems. The
deflation solves begin with an initial guess taken from a list of
(possibly several) initialization vectors. Newton iterations are
performed until a convergence tolerance is reached for a new solution
(added to the \emph{solution list} in Algorithm \ref{algo}) or a
maximum number of Newton iterations have been performed. For both the
deflated and undeflated Newton iterations, the convergence stopping
criterion on a given level is based on a set tolerance for an approximation's conformance to the first-order optimality conditions in the standard Euclidean $l^2$-norm. Throughout the numerical results section below, this tolerance is held at $10^{-4}$. For each Newton iteration, the linear systems are solved using the multigrid-preconditioned GMRES iterative solver proposed in \cite{Emerson3}, and the matrix-free approach outlined in Section \ref{DeflationMatrixFree} is applied when performing deflated iterations. Finally, the known solution approximations are transferred to a finer grid. In the current implementation, these finer grids represent successive uniform refinements of the initial coarse grid. 

\begin{algorithm}[h!]
\SetAlgoLined
~\\
0. Initialize $(\director_0, \phi_0, \lambda_0)$ on coarse grid.
~\\
1. Add guess, $(\director_0, \phi_0, \lambda_0)$, to continuation list.
~\\
\While{Refinement limit not reached}
{
	\tcp{First perform iterations on known solutions}
	\For{Each solution in continuation list}
	{
		\While{First-order optimality threshold not met}
		{
			2. Set up discrete \textbf{undeflated} system on current grid, $H$. ~\\
			3. Solve for $\ddirector_{H}$, $\dphi_{H}$, and $\dlambda_{H}$. ~\\
			4. Compute $\director_{k+1}$, $\phi_{k+1}$, and $\lambda_{k+1}$. ~\\
		}
		5. Add iterated solution to known solution list.
	}
	\tcp{Next perform deflation iterations to discover new solutions}
	6. Construct a set of initial guesses for the deflation solves on grid $H$. ~\\
	\For{Each guess in the list}
	{
		\While{First-order optimality threshold or failure criterion not met}
		{
			7. Set up discrete \textbf{deflated} system using known solutions on $H$. ~\\
			8. Solve for $\ddirector_{H}$, $\dphi_{H}$, and $\dlambda_{H}$. ~\\
			9. Compute $\director_{k+1}$, $\phi_{k+1}$, and $\lambda_{k+1}$. ~\\
		}
		\If{Conformance threshold met}
		{
			10. Add solution to known solution list. ~\\
		}
	}
	11. Uniformly refine the grid. ~\\
	12. Interpolate known solutions to fine grid $h$ and add to continuation list.
}
\caption{Newton's method minimization algorithm with NI and deflation}
\label{algo}
\end{algorithm}

For preconditioning of the linear solves, we use a geometric multigrid
implementation with standard finite-element interpolation operators
and Galerkin coarsening. This approach is monolithic, meaning that
coarse-grid correction treats all variables in the coupled system at
once; however, we use an efficient Braess-Sarazin-type relaxation
scheme that can be viewed as an approximate block factorization. A
single pre- and post-relaxation sweep is applied as part of a standard
V-cycle. Computational work for a full NI solve is given
in terms of work units (WUs), calculated as a weighted sum of the
total number of V-cycles across each NI level. With
uniform mesh refinements and a geometric multigrid strategy, the total
number of V-cycles on each grid is weighted by $(1/4)^l$, where $l$ is
the level of coarsening away from the finest mesh. For instance, the
total number of V-cycles on the second finest mesh is simply scaled by
$1/4$.  Thus, the total WUs for a given NI
solve provides a work measurement equivalent to counting fine-grid
V-cycles in a single-grid approach.

The blending of NI and deflation outlined above has a number of advantages above and beyond efficiency. Certain solutions are more readily detectable through a deflation process on a finer mesh. This is, for example, observed when considering cholesteric simulations or nematic configurations with defects for certain deflation parameters. Moreover, the algorithm allows for varying and adaptive initial guesses. That is, in addition to a static set of initial guesses for the deflation solves on each grid, sets of initial guesses may be constructed from transformations of known or newly discovered solutions throughout the NI and deflation process. Constructing strategies for dynamic generation of initial guesses will be the subject of future work.

\section{Numerical Results for Nematic Liquid Crystals} \label{NumericalExperiments}

In this section and the next, four numerical experiments using the
deflation approach detailed in Section \ref{DeflationMethod} are
carried out to demonstrate the performance of the method. The first
two simulations consider problems with known analytical solutions,
enabling validation of computed structures against the true
configurations. The remaining experiments illustrate the full
capabilities of the algorithm. For each simulation, the length scale
discussed in Section \ref{Model} is taken to be one micron, such that
$\sigma = 10^{-6}$ m. Furthermore, 
the characteristic Frank constant is taken to be $K = 6.2 \times
10^{-12}$ N, the dimensional value of $K_1$ for $5$CB, a common liquid
crystal, for convenience in adjusting relative parameter sizes. The applied
nondimensionalization, for instance, yields parameters $K_1 = 1$, $K_2 = 0.62903$,
and $K_3 = 1.32258$ for $5$CB. In addition, the characteristic voltage
is $\phi_0 = 1$ V, which implies that the nondimensional dielectric permittivity constant is $\epsilon_0 = 1.42809$.

For the test problems, we consider a classical domain with two
parallel substrates placed at unit distance apart. These substrates
run parallel to the $xz$-plane and perpendicular to the
$y$-axis. Further, we assume a slab-type domain such that $\director$
may have a non-zero $z$-component, but $\pd{\director}{z} =
\vec{0}$. Thus, for the numerical experiments to follow, $\Omega = \{
(x,y) \text{ } : \text{ } 0 \leq x,y \leq 1 \}$.  For the first
two experiments, periodic boundary conditions are applied at the left
and right boundaries and Dirichlet conditions are enforced at the top
and bottom of the domain. In the third experiment, Dirichlet boundary
conditions are applied for the entire boundary. The deflation
parameters are fixed such that $\alpha = 1$ and $p = 3$, and the
failure criterion in Algorithm \ref{algo} occurs when the number of
Newton iterations reaches $100$ or the average length of the current
director field is above $3.0$, substantially violating the unit-length
constraint.  (Similar failure criteria could be integrated with the
first Newton loop in Algorithm 1, for continuing known solutions, but
this appears to be unnecessary in practice.)  In each simulation, the algorithm begins on a uniform $8 \times 8$ coarse mesh, ascending in uniform refinements to a $256 \times 256$ fine grid. The algorithm's discretizations and grid management are performed with the deal.II scientific computing library \cite{BangerthHartmannKanschat2007}. In the numerical tests to follow, biquadratic finite elements are used to discretize components associated with $\director$ and $\phi$ in the variational systems above, while piecewise constants are used for those related to $\lambda$. This results in a total of $1,118,212$ degrees of freedom on the finest mesh.

A form of damped Newton stepping is applied for both the undeflated and deflated iterate updates such that the updated iterates are given by
\begin{align*}
\left[ \begin{array}{c} 
\director_{k+1} \\
\phi_{k+1} \\
\lambda_{k+1}
\end{array} \right] =
\left[ \begin{array}{c} 
\kdirector \\
\kphi \\
\klambda
\end{array} \right] + \omega
\left [ \begin{array}{c}
\ddirector \\
\dphi \\
\dlambda
\end{array}\right],
\end{align*}
where $0 < \omega \leq 1$ is a damping parameter. For the undeflated
solves, $\omega = \omega_1$ on the coarse grid and is increased by
$\Delta_1$ at each refinement to a maximum of $1$. With the deflated
systems, $\omega = \omega_2$ on the coarse grid and decreases by
$\Delta_2$ to a minimum of $0.1$. This strategy aims at improving
convergence for both types of iterations. The damping parameter is
increased on each grid for the undeflated solves as confidence in the
Newton convergence increases for more finely resolved solutions. On
the other hand, at each level, all of the deflation iterations begin
with a relatively na\"{i}ve initial guess and the possibility of
convergence becomes more tenuous on finer meshes. Hence, the
decreasing damping parameter invokes tighter control over step length
on finer grids. Note that more advanced step selection methods, such
as trust regions \cite{Emerson3, Byrd2}, exist. However,
we experimentally observed that using trust regions during the
deflation phase of the algorithm hindered the methods ability to
discover new basins of attraction, thereby limiting the number of
unique solutions found. Improving this performance will be considered
in future work. Finally, the linear solver tolerance, which is
based on a ratio of the norm of the current (discrete) solution's residual to that of the initial guess, is held at $10^{-6}$.

\subsection{Tilt-Twist Configuration} \label{TiltTwistExperiments}

The first problem considered in this section is an elastic configuration with no electric field and Frank constants given by $K_1 = 1.0$, $K_2 = 3.0$, and $K_3 = 1.2$. For the Newton damping, $\omega_1 = 1.0$, $\Delta_1 = 0.0$, $\omega_2 = 1.0$, and $\Delta_2 = 0.5$. At the Dirichlet boundaries, we set
\begin{align*}
\director(x, 0) &= \left( \cos \left(-\frac{\pi}{4} \right), 0, \sin \left(-\frac{\pi}{4}\right) \right), & & \director(x, 1) = \left( \cos \left(\frac{\pi}{4}\right), 0, \sin \left(\frac{\pi}{4}\right) \right).
\end{align*}
This is known as a tilt-twist problem and is an interesting example for a few reasons. The opposing boundary conditions induce a twisting configuration in the nematics through the interior of the domain. Under these conditions, a planar twisting pattern, where the $y$-component of the director remains zero, satisfies the first-order optimality conditions. However, for these Frank constants it is well known that a twist configuration incorporating a nonplanar tilt is energetically optimal \cite{Stewart1, Leslie2}. Thus, there are multiple solutions satisfying \eqref{FOOC1}-\eqref{FOOC3}. Furthermore, these nonplanar twist solutions only become energetically optimal for certain Frank constant ratios. For instance, such configurations are not detectable when using the one-constant approximation \cite{Stewart1, Cohen1}.

For the deflation solves, two initial guesses are constructed at each refinement level to serve as starting points for the discovery of additional solutions. Through the interior of the domain, both initial guesses are isolated to the $xy$-plane and incorporate a slight uniform tilt; see Appendix \ref{InitialGuesses}. As discussed in \cite{Emerson3}, convergence to the energetically optimal solution can be attained even when choosing a relatively na\"{i}ve initial guess. However, without deflation, the poor initial guesses used here result in convergence on all grids to a single planar twist solution, which represents only a local minimum. The first guess is also used for the coarse-grid, undeflated iterations.

\begin{figure}[h!]
\centering
\begin{subfigure}{.32 \textwidth}
  \includegraphics[scale=.205, left]{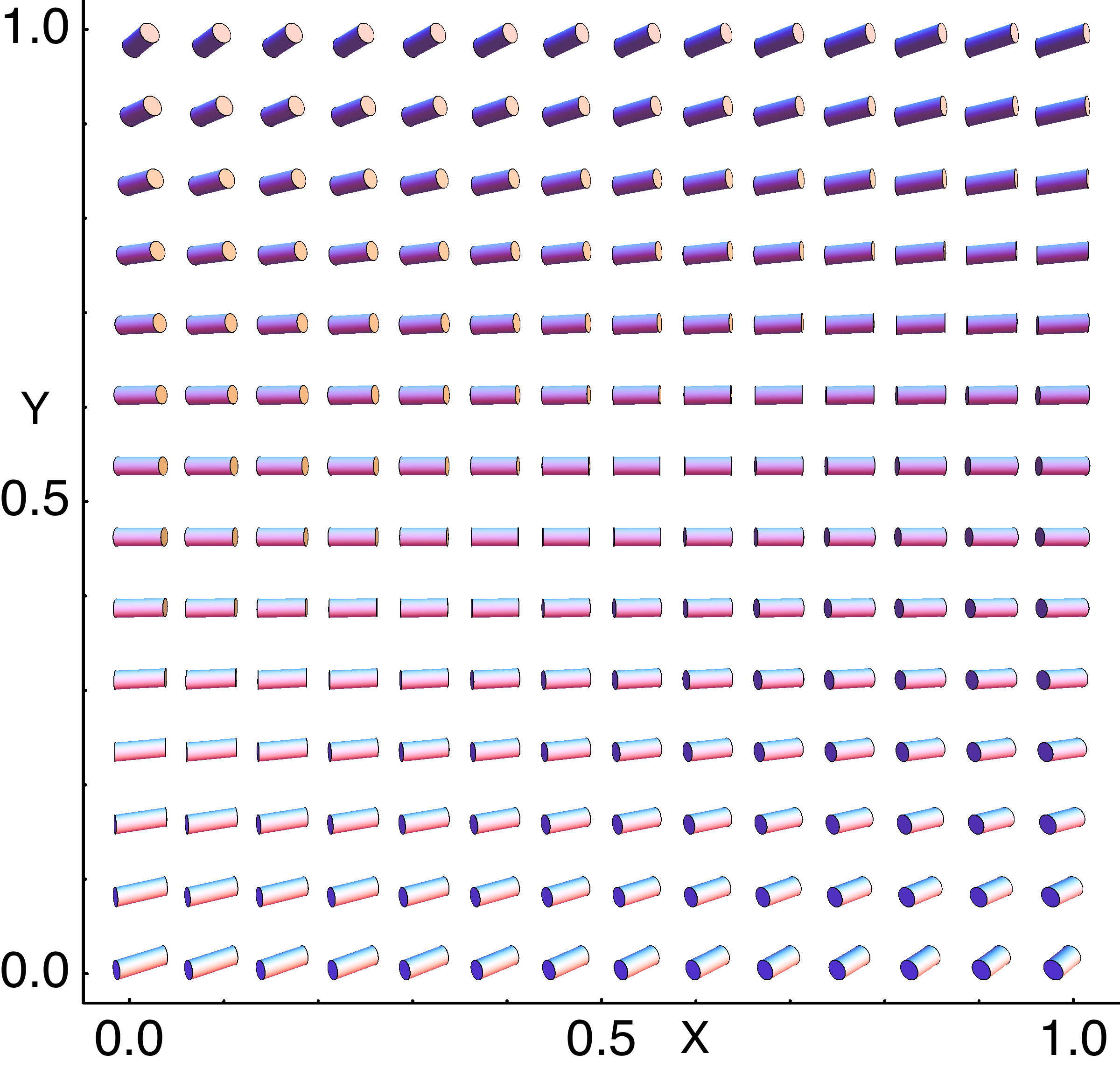}
    \caption{}
  \label{TiltTwistSolutions:left1}
\end{subfigure}
\begin{subfigure}{.32 \textwidth}
  \includegraphics[scale=.205, center]{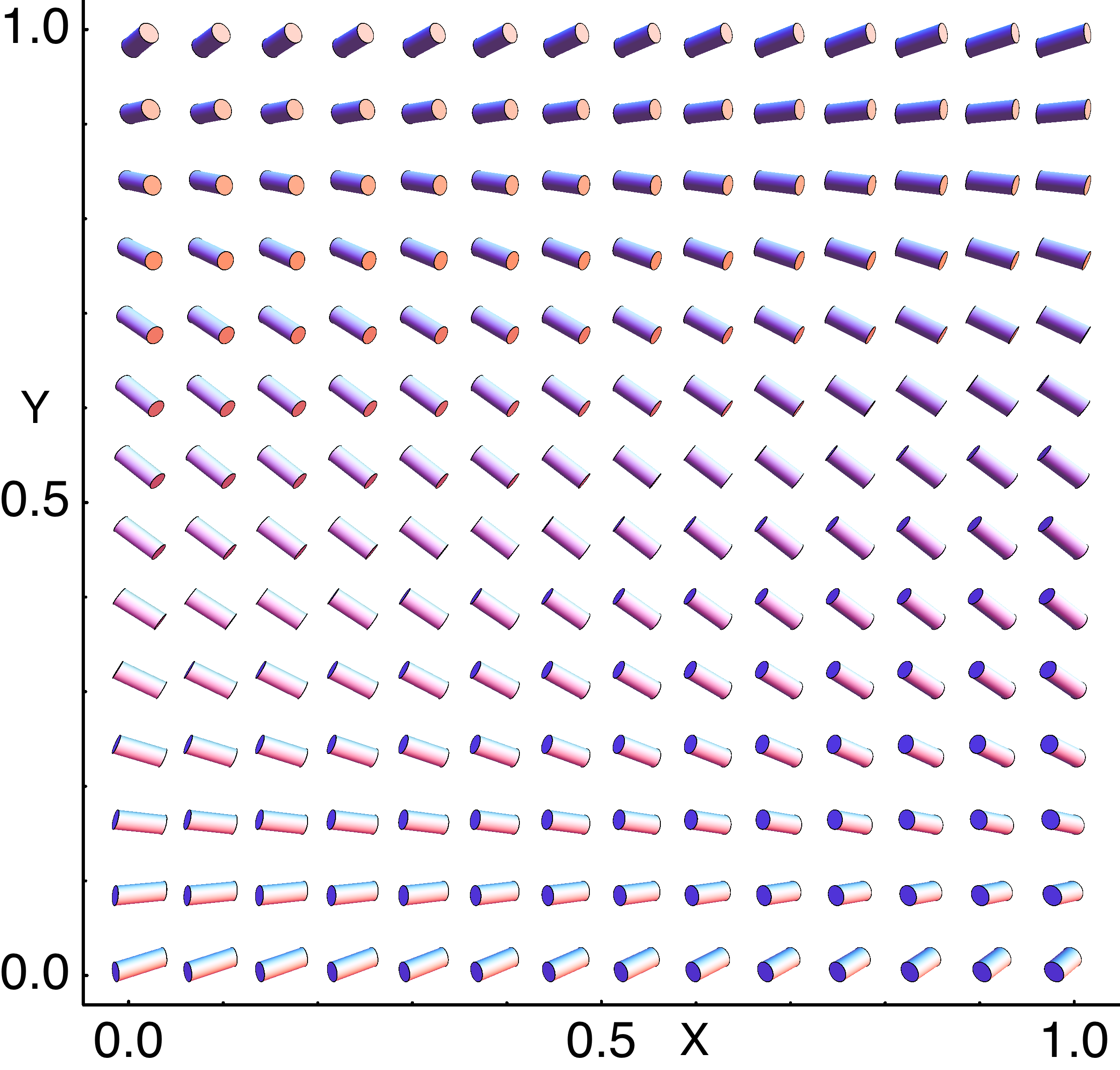}
      \caption{}
  \label{TiltTwistSolutions:center1}
\end{subfigure}
\begin{subfigure}{.32 \textwidth}
        \includegraphics[scale=.205, right]{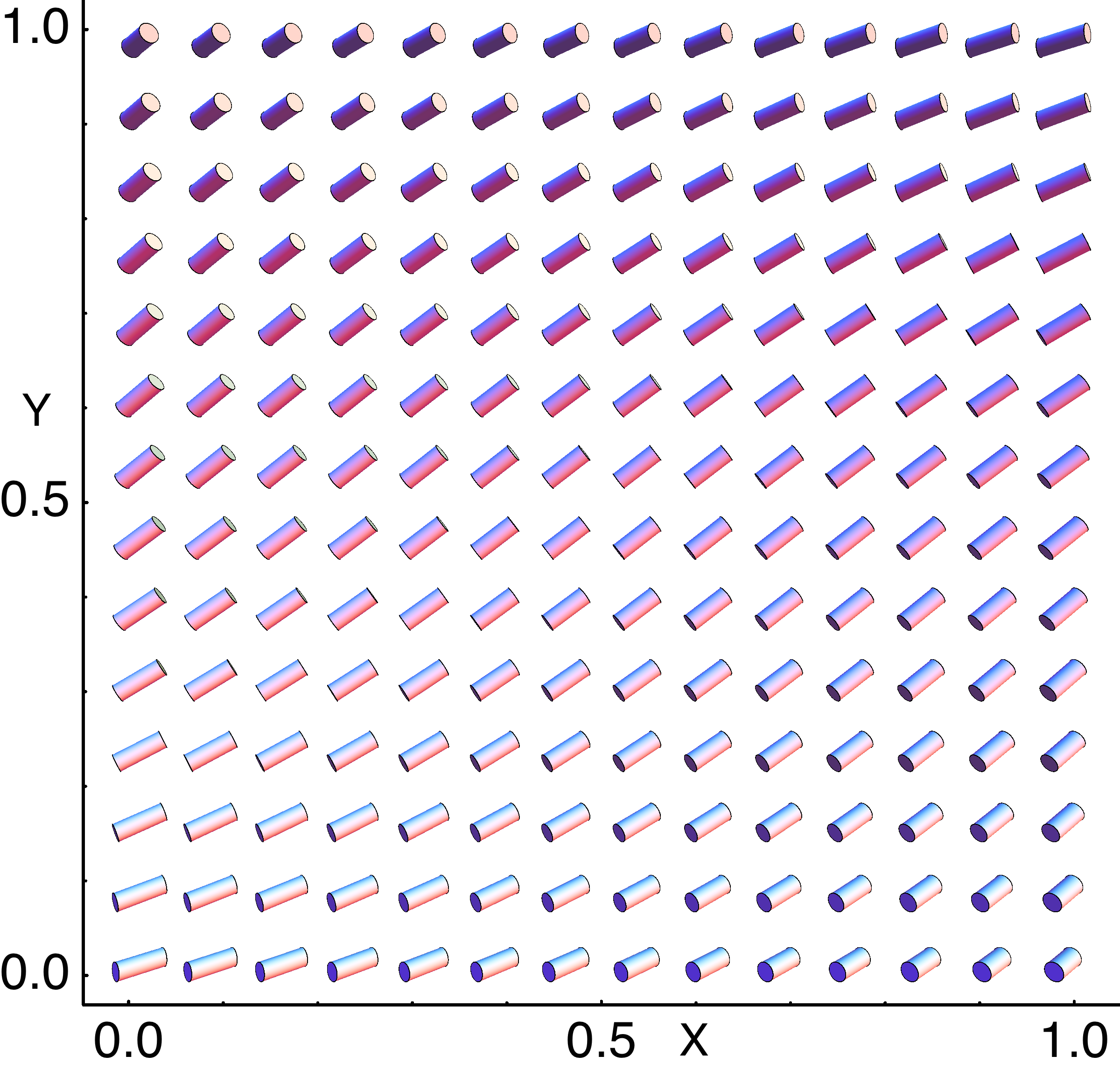}
      \caption{}
  \label{TiltTwistSolutions:right1}
\end{subfigure}
\begin{subfigure}{.32 \textwidth}
  \includegraphics[scale=.205, center]{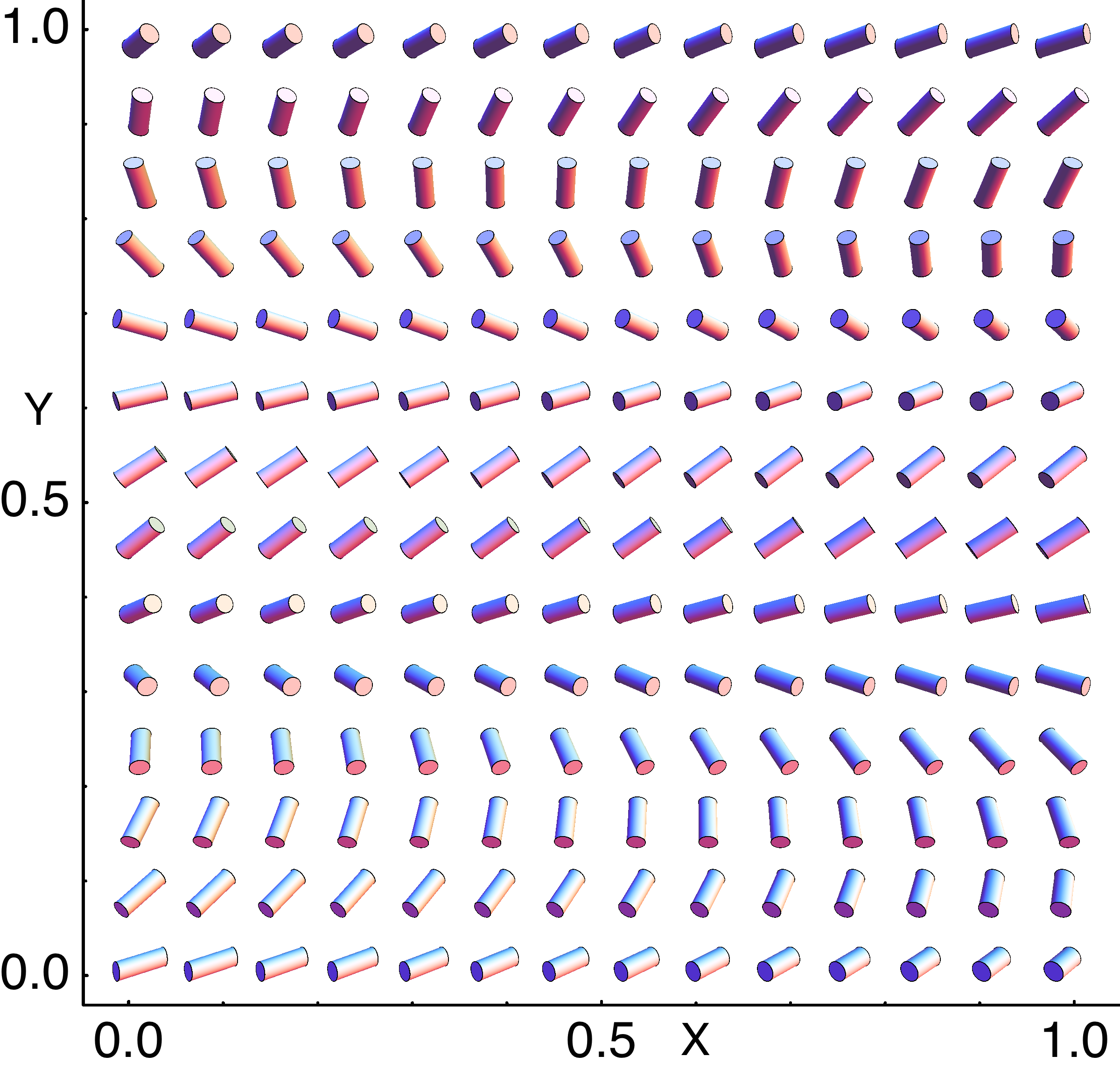}
      \caption{}
  \label{TiltTwistSolutions:center2}
\end{subfigure}
\caption{\small{(\subref{TiltTwistSolutions:left1}) Resolved non-minimizing solution on $256 \times 256$ mesh (restricted for visualization) with final free energy of $3.701$. (\subref{TiltTwistSolutions:center1}) Energy-minimizing solution identified through deflation with final free energy of $3.593$. (\subref{TiltTwistSolutions:right1}) Symmetric energy-minimizing solution found with deflation. All solutions were located on the coarsest mesh. (\subref{TiltTwistSolutions:center2}) Intricate non-minimizing solution satisfying the optimality conditions located with deflation parameters $\alpha = 0.1$ and $p = 2.0$.}}
\label{TiltTwistSolutions}
\end{figure}

The undeflated iterations converge to the planar twist solution displayed in Figure \ref{TiltTwistSolutions}(\subref{TiltTwistSolutions:left1}) with a final free energy of $3.701$ and consume a total of $11.9$ WUs as the solution is continued through the NI hierarchy. The remainder of the solutions are located using deflation. The configurations in Figure \ref{TiltTwistSolutions}(\subref{TiltTwistSolutions:center1}) and (\subref{TiltTwistSolutions:right1}) represent the energetically optimal structures for this problem, with both exhibiting final free energies of $3.593$. Newton solves for these solutions require $17.8$ and $18.0$ WUs, respectively. The symmetry of the device and boundary conditions allow for the reflection symmetry seen in these two solutions. As discussed above, certain selections of deflation parameters yield additional distinct solutions. The configuration displayed in Figure \ref{TiltTwistSolutions}(\subref{TiltTwistSolutions:center2}) is discovered when applying deflation parameters of $\alpha = 0.1$ and $p = 2.0$. The associated free energy is $32.336$. While the configuration is clearly not energetically optimal, it satisfies the first-order optimality conditions. 

\begin{figure}[h!]
\begin{subfigure}{.49 \textwidth}
  \includegraphics[scale=.32, left]{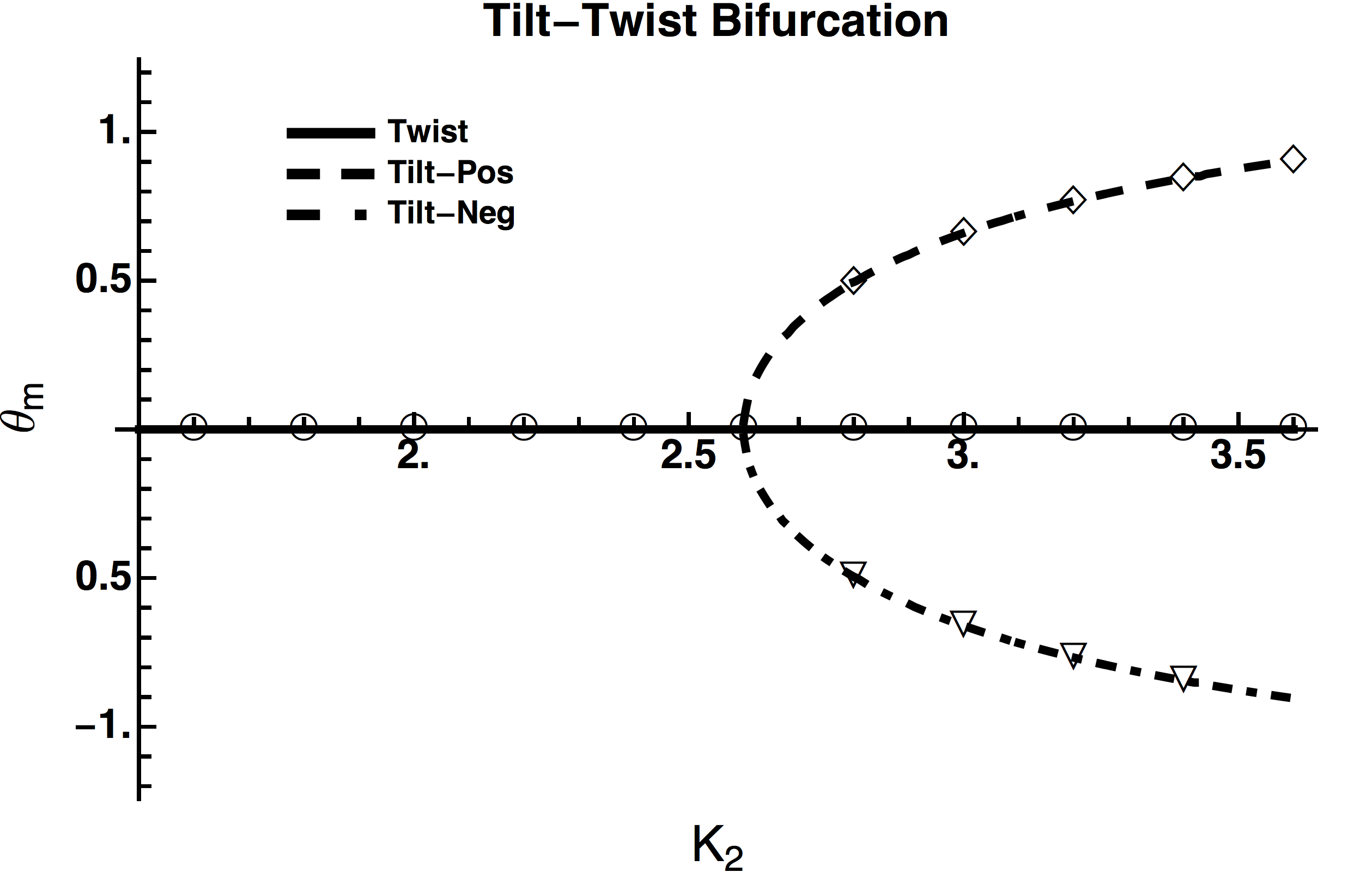}
    \caption{}
  \label{TiltTwistBifurcation:left1}
\end{subfigure} \hfill
\begin{subfigure}{.49 \textwidth}
  \includegraphics[scale=.32, right]{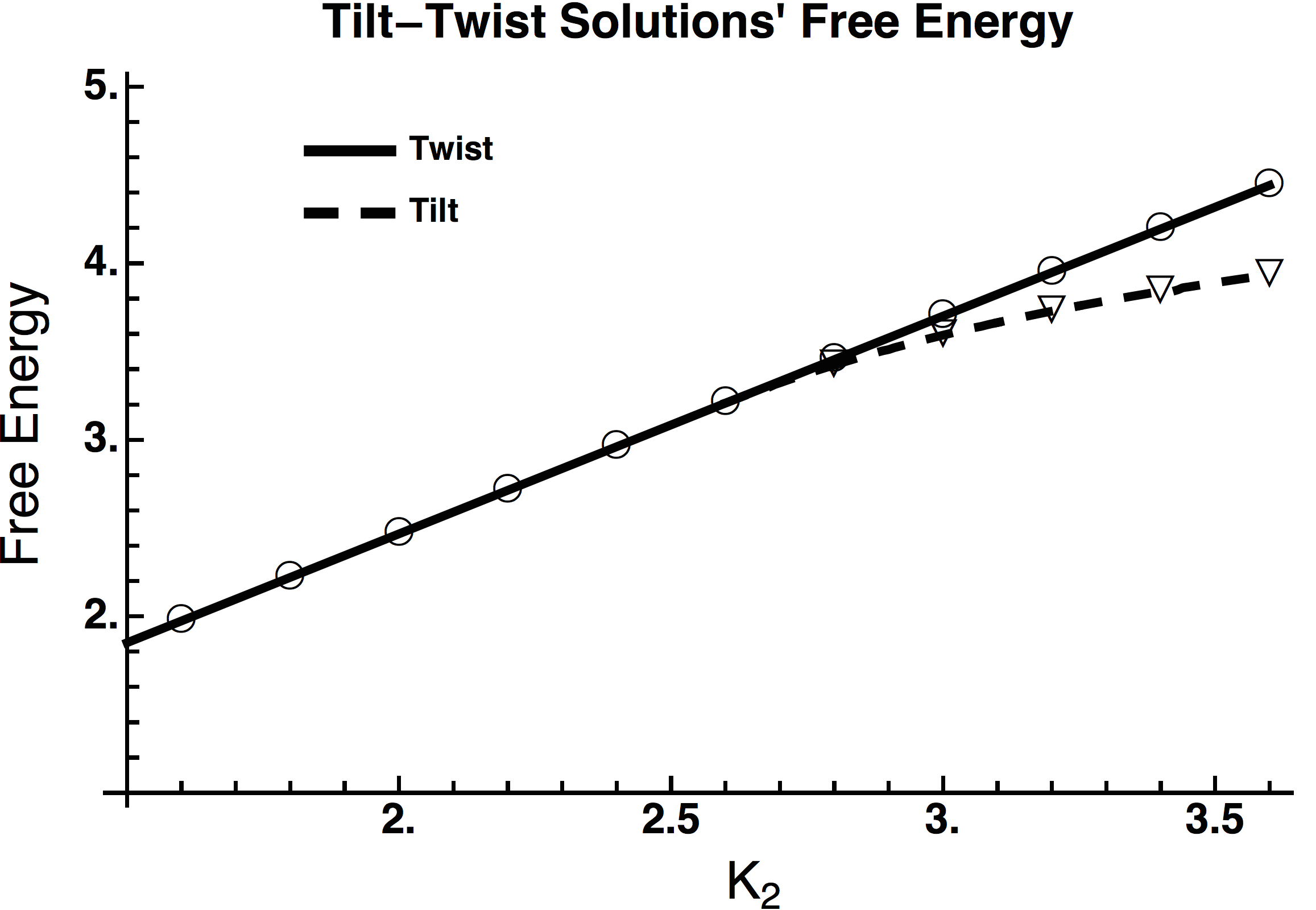}
    \caption{}
  \label{TiltTwistBifurcation:right1}
\end{subfigure}
\caption{\small{(\subref{TiltTwistBifurcation:left1}) Pitchfork bifurcation diagram for the tilt-twist problem with varying $K_2$ values. Lines depict analytical values for $\theta_m$ while markers indicate maximum angular tilt for solutions obtained through deflation. (\subref{TiltTwistBifurcation:right1}) Free energy as a function of $K_2$. Lines are analytical free energies and markers denote free energies for solutions obtained through deflation.}}
\label{TiltTwistBifurcationDiagrams}
\end{figure}

The behavior of the solution branches for this problem can be characterized by the value of the Frank constant, $K_2$, if the remainder of the parameters are held fixed. As the magnitude of $K_2$ varies, the admissible solution set undergoes a pitchfork bifurcation process. This bifurcation delineates the transition away from energetic optimality of the simple twist solution towards tilt-twist solutions like that seen in Figure \ref{TiltTwistSolutions}(\subref{TiltTwistSolutions:center1}). This behavior is captured in Figure \ref{TiltTwistBifurcationDiagrams}(\subref{TiltTwistBifurcation:left1}). Holding the boundary conditions and remaining Frank values constant, the figure displays the maximal azimuthal angle, $\theta_m$, in the liquid crystal configuration as a function of $K_2$. This deviation describes the extent of ``tilt" present in a given solution. The lines represent the known analytical value of $\theta_m$ \cite{Stewart1, Leslie2} for the different solution branches, while the individual markers correspond to values of $\theta_m$ computed for each of the solutions located with the deflation algorithm.

Increasing the value of $K_2$ increases the energetic contribution of twist-type deformations in the liquid crystal structure. Thus, as $K_2$ rises, so does the incentive to reduce twisting through tilt. When $K_2$ reaches a critical value, this incentive is large enough to produce an energetically optimal solution with nonzero tilt. For the set of constants considered, this occurs at approximately $K_2 = 2.61$, after which at least three solutions satisfy the first-order optimality conditions. The first is the simple twist solution, while the second and third are energetically optimal tilt-twist solutions related by reflection. Figure \ref{TiltTwistBifurcationDiagrams}(\subref{TiltTwistBifurcation:right1}) displays the energetic behavior of the simple twist and tilt-twist solutions for increasing values of $K_2$. There, the point of separation for the free energies is clearly visible, with the minimal free energy transitioning to the tilt-twist arrangements above the critical threshold. The lines represent the known, analytical free energy and individual markers are computed free energies using the deflation process.

\vspace{0.1cm}

\subsection{Freedericksz Transition} \label{FreederickszExperiments}

The second numerical experiment considers a classical Freedericksz transition problem with simple director boundary conditions such that $\director$ lies uniformly parallel to the $x$-axis at the edges $y=0$ and $y=1$. For the electric potential, $\phi$, the boundary conditions set $\phi(x,0) = 0$ and $\phi(x,1) = V = 1.1$. The relevant Frank and electric constants are $K_1 = 1$, $K_2=0.62903$, and $K_3= 1.32258$, $\epsilon_0 = 1.42809$, $\epsilon_{\perp}=7$, and $\epsilon_a = 11.5$. Note that for $\epsilon_a > 0$ the liquid crystals are attracted to alignment parallel to the electric field. The relevant damping parameters are $\omega_1 = 1.0$, $\Delta_1 = 0.0$, $\omega_2 = 1.0$, and $\Delta_2 = 0.5$. The same two initial guesses for $\director$ used in the previous experiment are applied in the deflation solves here; c.f. Appendix \ref{InitialGuesses}. These configurations serve as the starting point for all deflation searches in the NI hierarchy.

\begin{figure}[h!]
\centering
\begin{subfigure}{.32 \textwidth}
  \includegraphics[scale=.205, left]{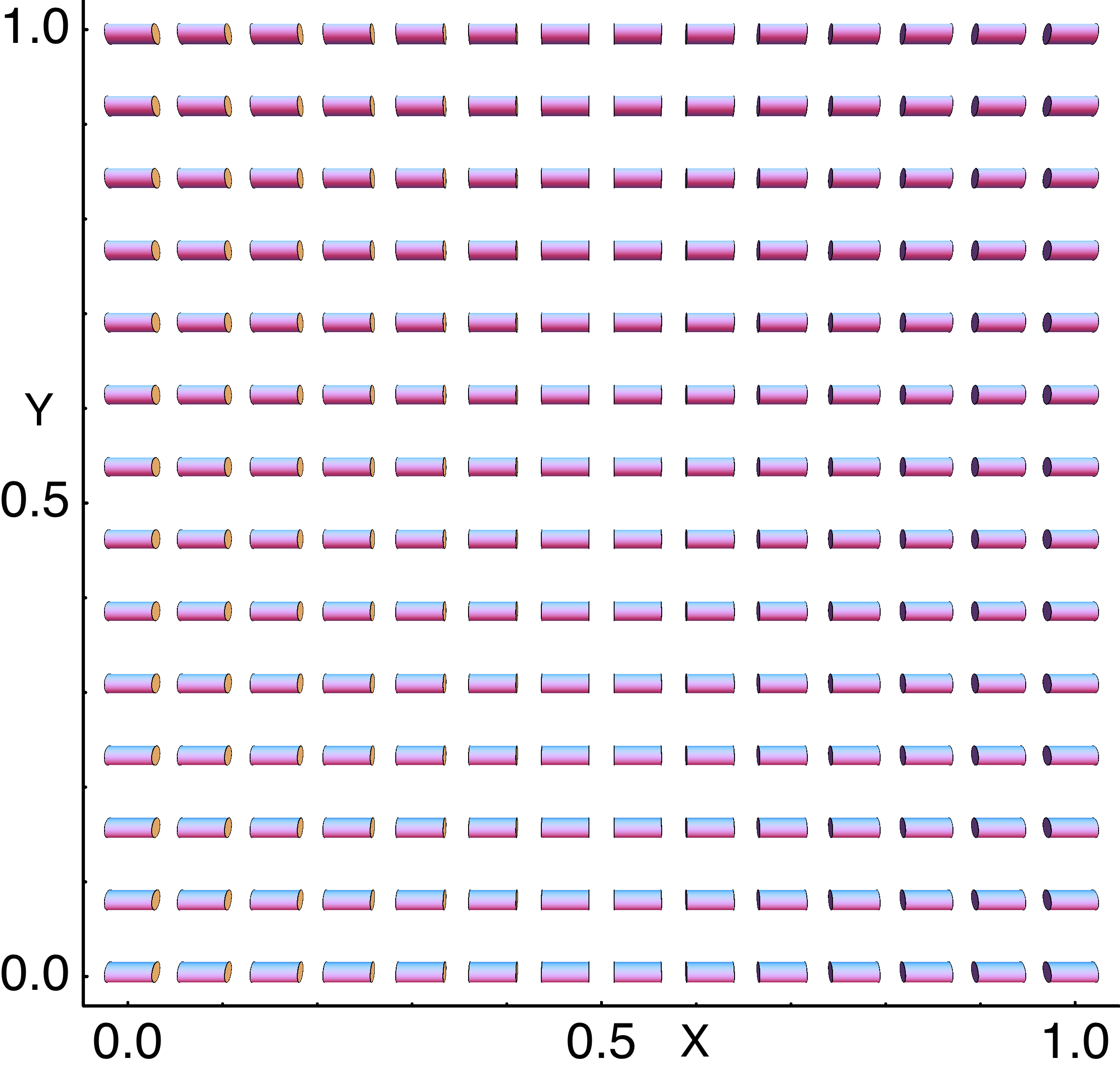}
    \caption{}
  \label{FreederickszTransSolutions:left1}
\end{subfigure}
\begin{subfigure}{.32 \textwidth}
  \includegraphics[scale=.205, center]{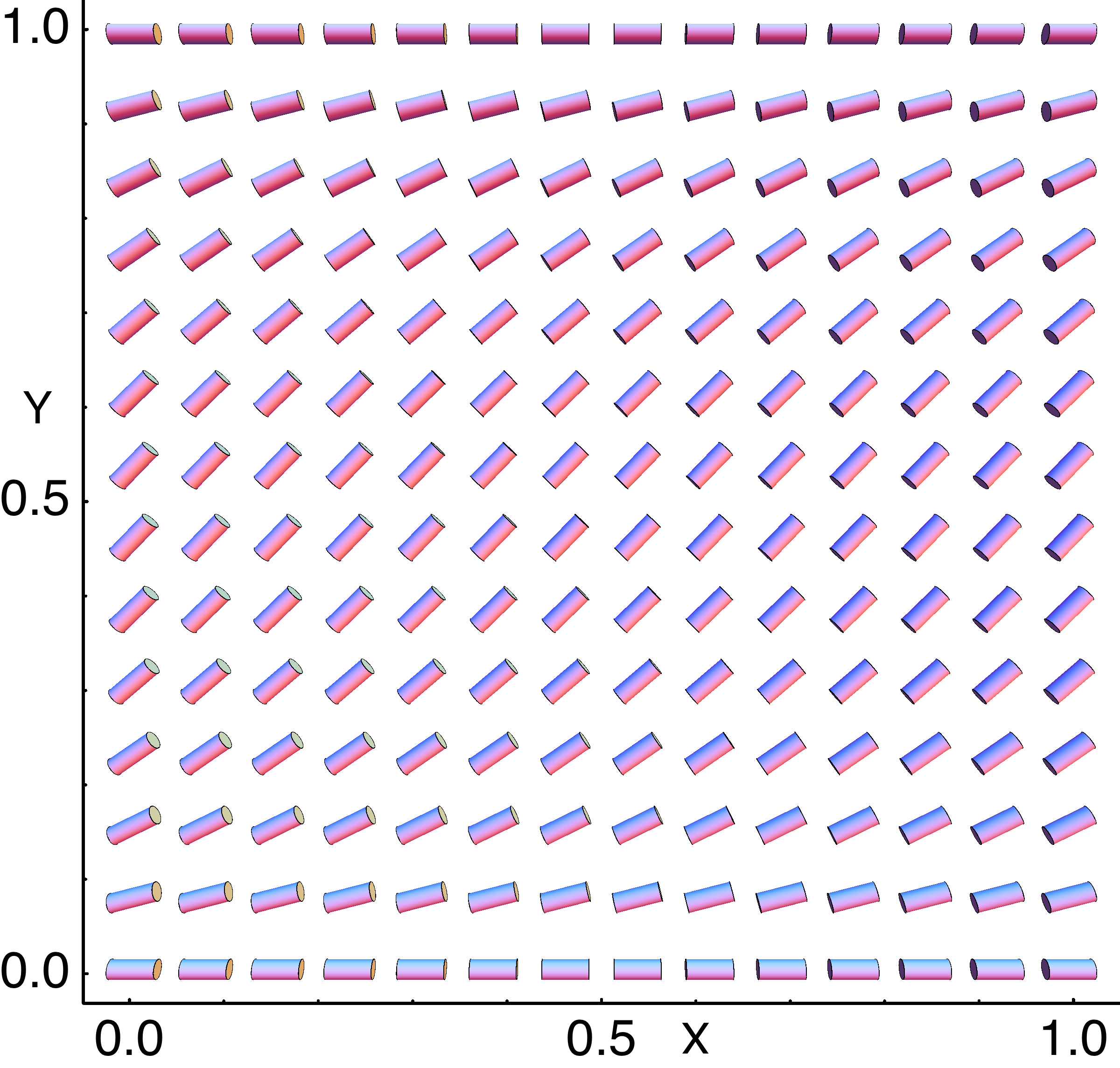}
      \caption{}
  \label{FreederickszTransSolutions:center1}
\end{subfigure}
\begin{subfigure}{.32 \textwidth}
        \includegraphics[scale=.205, right]{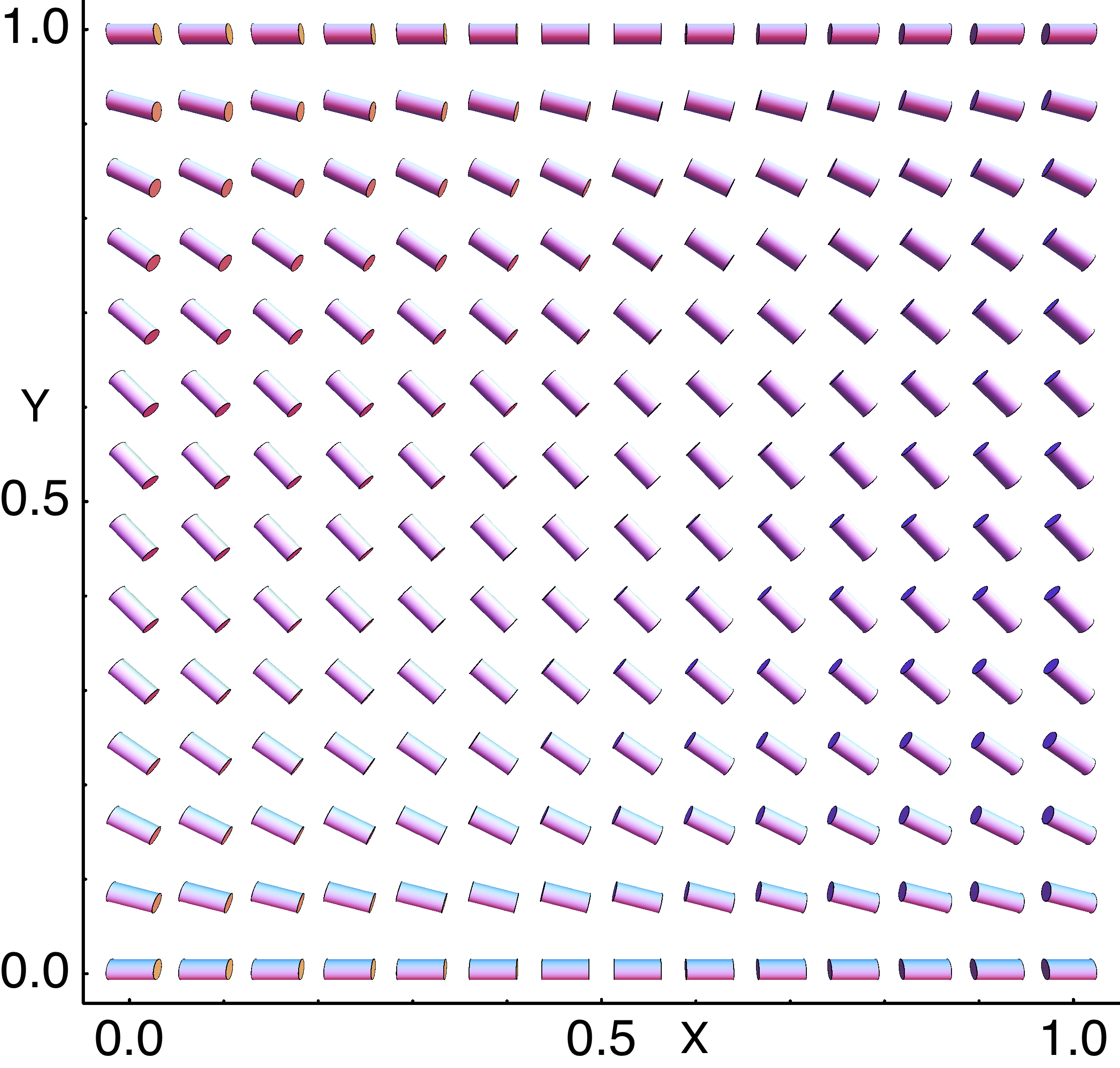}
      \caption{}
  \label{FreederickszTransSolutions:right1}
\end{subfigure}
\caption{\small{(\subref{FreederickszTransSolutions:left1}) Resolved non-minimizing solution on $256 \times 256$ mesh (restricted for visualization) with final free energy of $-6.048$. (\subref{FreederickszTransSolutions:center1}) Energy-minimizing solution found through deflation with final free energy of $-6.778$. (\subref{FreederickszTransSolutions:right1}) Symmetric energy-minimizing solution computed with deflation. All solutions were located on the coarsest mesh.}}
\label{FreederickszTransSolutions}
\end{figure}

The initial undeflated iterations converge to the elastic rest configuration uniformly parallel to the $x$-axis shown in Figure \ref{FreederickszTransSolutions}(\subref{FreederickszTransSolutions:left1}) and use $16.0$ WUs. The final free energy for this structure is $-6.048$. Thereafter, using deflation, the energetically optimal arrangements displayed in Figures \ref{FreederickszTransSolutions}(\subref{FreederickszTransSolutions:center1}) and (\subref{FreederickszTransSolutions:right1}) are found and each has a final free energy of $-6.778$. The computation of each solution requires $33.4$ WUs. These solutions represent a true Freedericksz transition in which the applied electric field successfully deforms the nematic configuration away from elastic rest. Without deflation the two guesses used here converge to the same solution, Figure \ref{FreederickszTransSolutions}(\subref{FreederickszTransSolutions:left1}), on each grid.

As with the tilt-twist configurations, the Freedericksz transition problem exhibits an important pitchfork bifurcation. The strength of the applied voltage at the top substrate, $V$, relative to the elastic characteristics of the liquid crystal, determines the bifurcation structure. Retaining the liquid crystal constants outlined above and varying the applied voltage, we observe the bifurcation process. As the applied voltage becomes stronger, the electric field begins to overpower the elastic effects in the sample. At a critical threshold, given analytically by $V_c = \pi \sqrt{\frac{K_1}{\epsilon_0 \epsilon_a}}$, it becomes energetically advantageous to tilt in the direction of the field \cite{Zocher1, Stewart1}. The critical voltage for the problem parameters considered here is $V_c = 0.775$.

\begin{figure}[h!]
\begin{subfigure}{.49 \textwidth}
  \includegraphics[scale=.315, left]{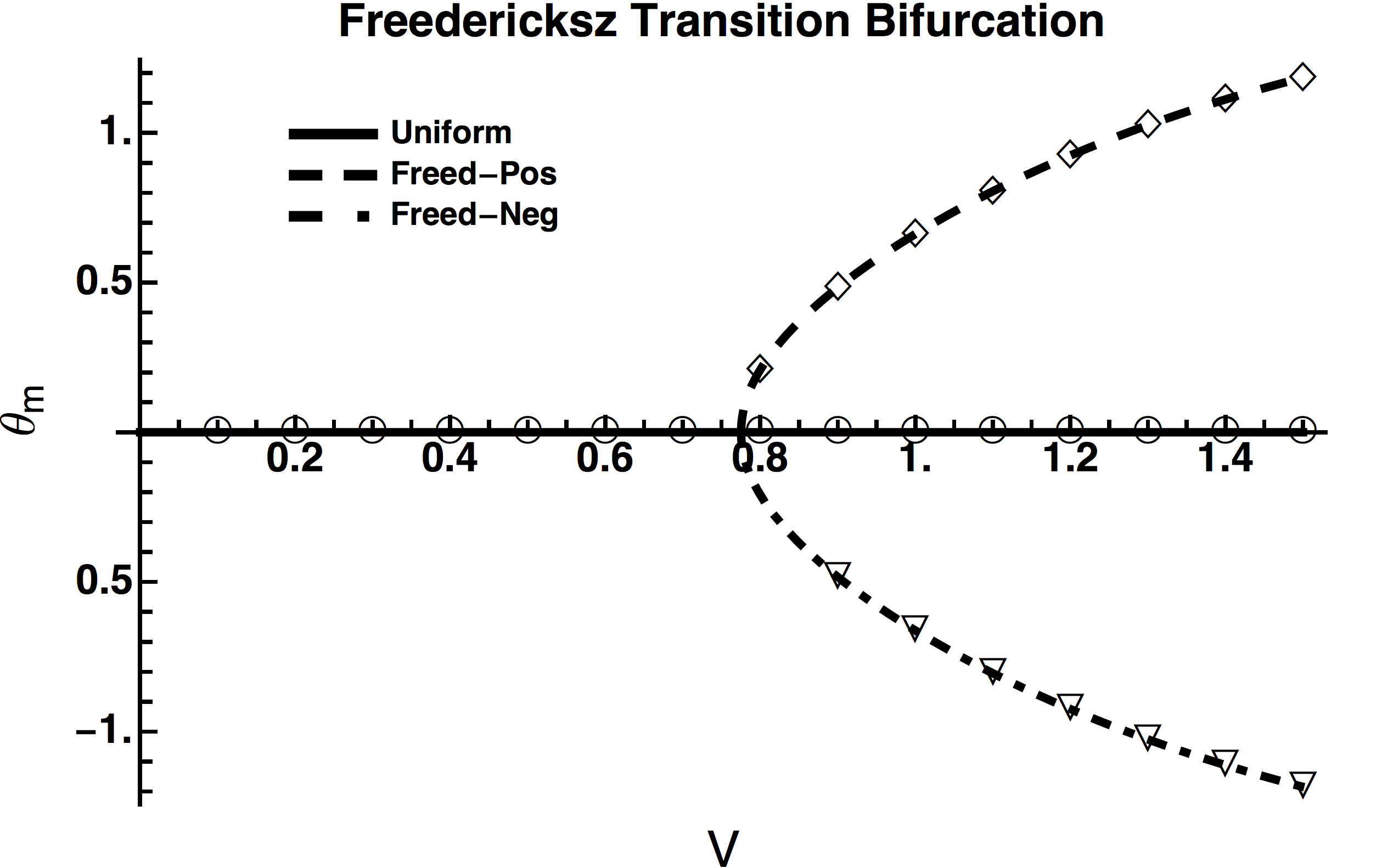}
    \caption{}
  \label{FreederickszBifurcation:left1}
\end{subfigure} \hfill
\begin{subfigure}{.49 \textwidth}
  \includegraphics[scale=.315, right]{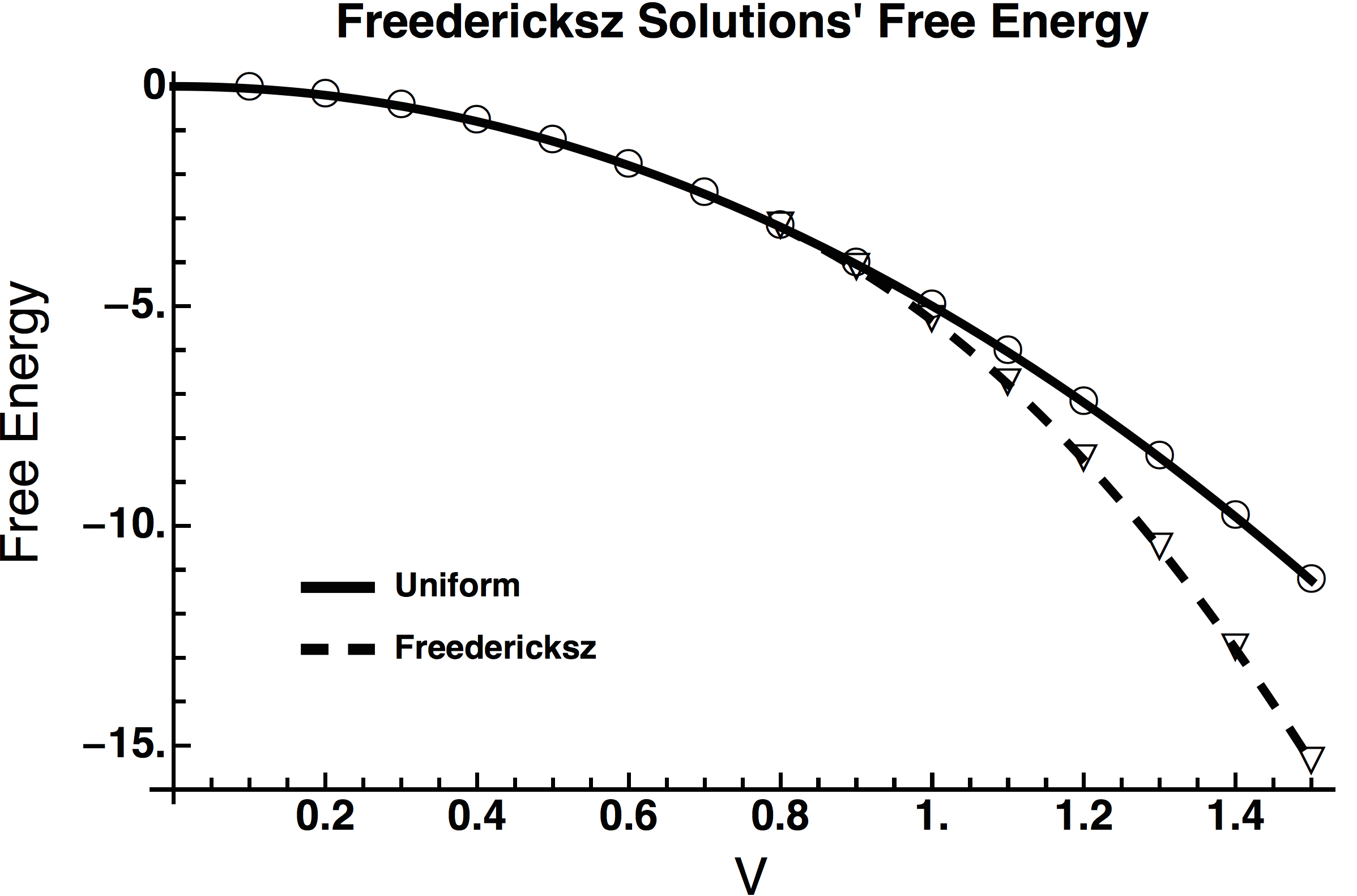}
    \caption{}
  \label{FreederickszBifurcation:right1}
\end{subfigure}
\caption{\small{(\subref{FreederickszBifurcation:left1}) Pitchfork bifurcation diagram characterizing the Freedericksz transition at approximately $V_c = 0.775$. Lines depict analytical values for $\theta_m$ while markers indicate maximum angular tilt for solutions obtained through deflation. (\subref{FreederickszBifurcation:right1}) A plot of free energy as a function of applied voltage. Lines are analytical free energies and markers denote free energies for solutions obtained through deflation.}}
\label{FreederickszBifurcationDiagrams}
\end{figure}

In Figure \ref{FreederickszBifurcationDiagrams}(\subref{FreederickszBifurcation:left1}), when $V$ reaches the critical value, solutions tilting in the direction of the electric field begin to satisfy the first-order optimality conditions and yield optimal free energy. The value $\theta_m$ denotes the maximum angular tilt of the director field in the direction of the electric field resulting from the applied voltage. As in Figure \ref{TiltTwistBifurcationDiagrams}(\subref{TiltTwistBifurcation:left1}), the lines represent analytical computations for $\theta_m$ as $V$ varies \cite{Stewart1}, and the individual markers are values for solutions computed through deflation. Figure \ref{FreederickszBifurcationDiagrams}(\subref{FreederickszBifurcation:right1}) characterizes the shift in free-energy optimizing solutions resulting from the Freedericksz transition as $V$ passes the critical voltage, $V_c$. The lines represent known, analytical free energies while the free energies associated with solutions computed through the deflation procedure are denoted with the individual markers. 

\subsection{Escape and Disclination Solutions} \label{EscapeDisclination}

This third numerical experiment investigates the phenomenon of defects, also known as disclinations. Defects in liquid crystal structures are locations in a sample where the director field is undefined or contains discontinuities. There are a multitude of disclination types including point, wedge, sheet, and loop defects, among others. In this example, we consider wedge disclinations. These disclinations involve rotation around an axis parallel to the defect and are, therefore, sometimes referred to as axial disclinations \cite{Friedel1}. Wedge-type disclinations have been studied in \cite{Frank1, Dzyaloshinskii1}. 

For this simulation, the damping parameters are $\omega_1 = 0.4$, $\Delta_1 = 0.2$, $\omega_2 = 1.0$, and $\Delta_2 = 0.5$. Dirichlet boundary conditions are applied to the entire domain boundary and no electric field is present. The boundary conditions are fixed such that the director faces the center of the domain. We use the Frank constants $K_1 = 1.0$, $K_2 = 3.0$, and $K_3 = 1.2$. As in the previous experiments, two initial guesses, detailed in Appendix \ref{InitialGuesses}, are used for the deflation solves on each grid.

\begin{figure}[h!]
\centering
\begin{subfigure}{.32 \textwidth}
  \includegraphics[scale=.205, left]{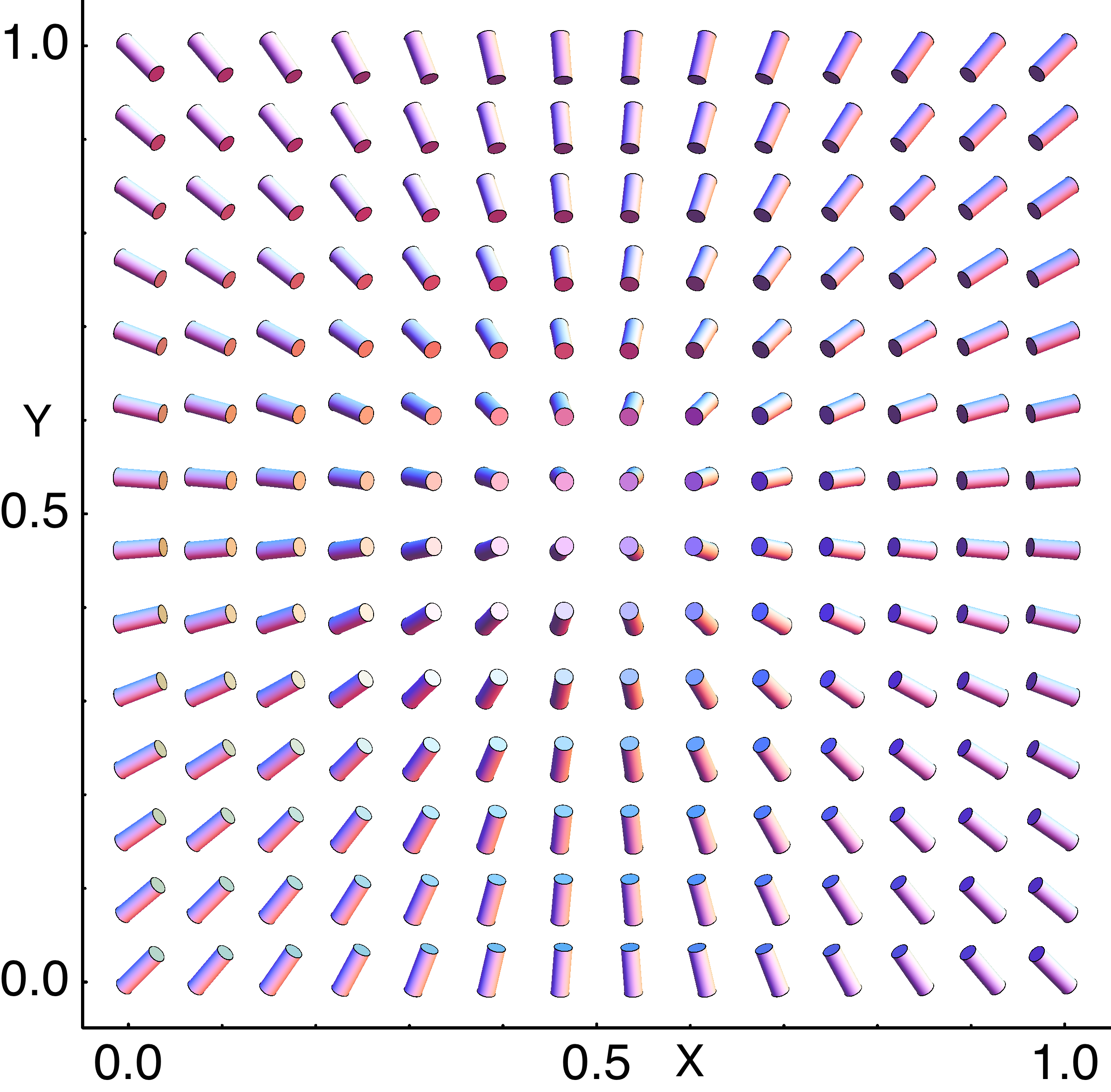}
    \caption{}
  \label{EscapeDisclinationConfig:left1}
\end{subfigure}
\begin{subfigure}{.32 \textwidth}
  \includegraphics[scale=.205, center]{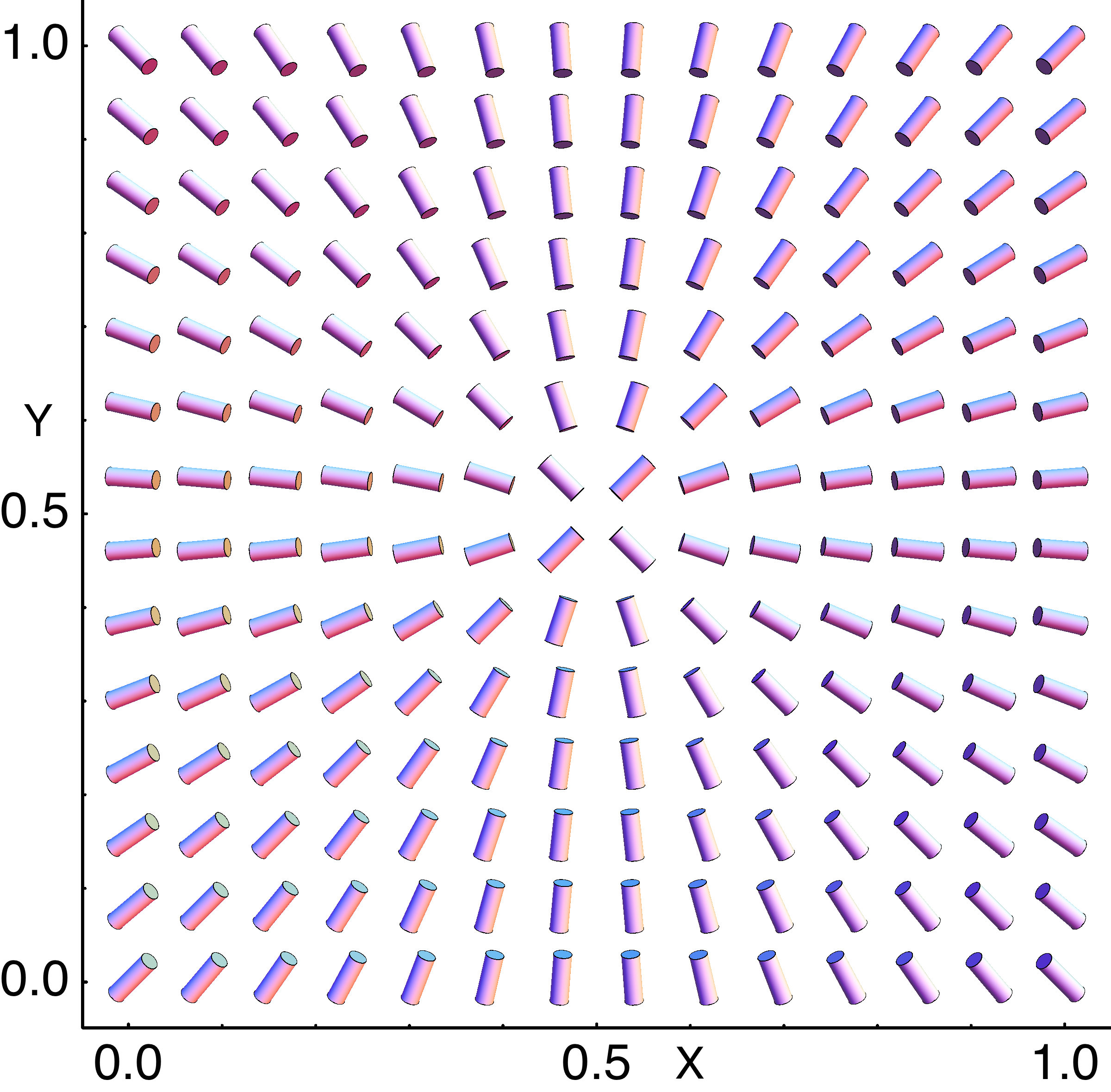} 
      \caption{}
  \label{EscapeDisclinationConfig:center1}
\end{subfigure}
\begin{subfigure}{.32 \textwidth}
        \includegraphics[scale=.205, right]{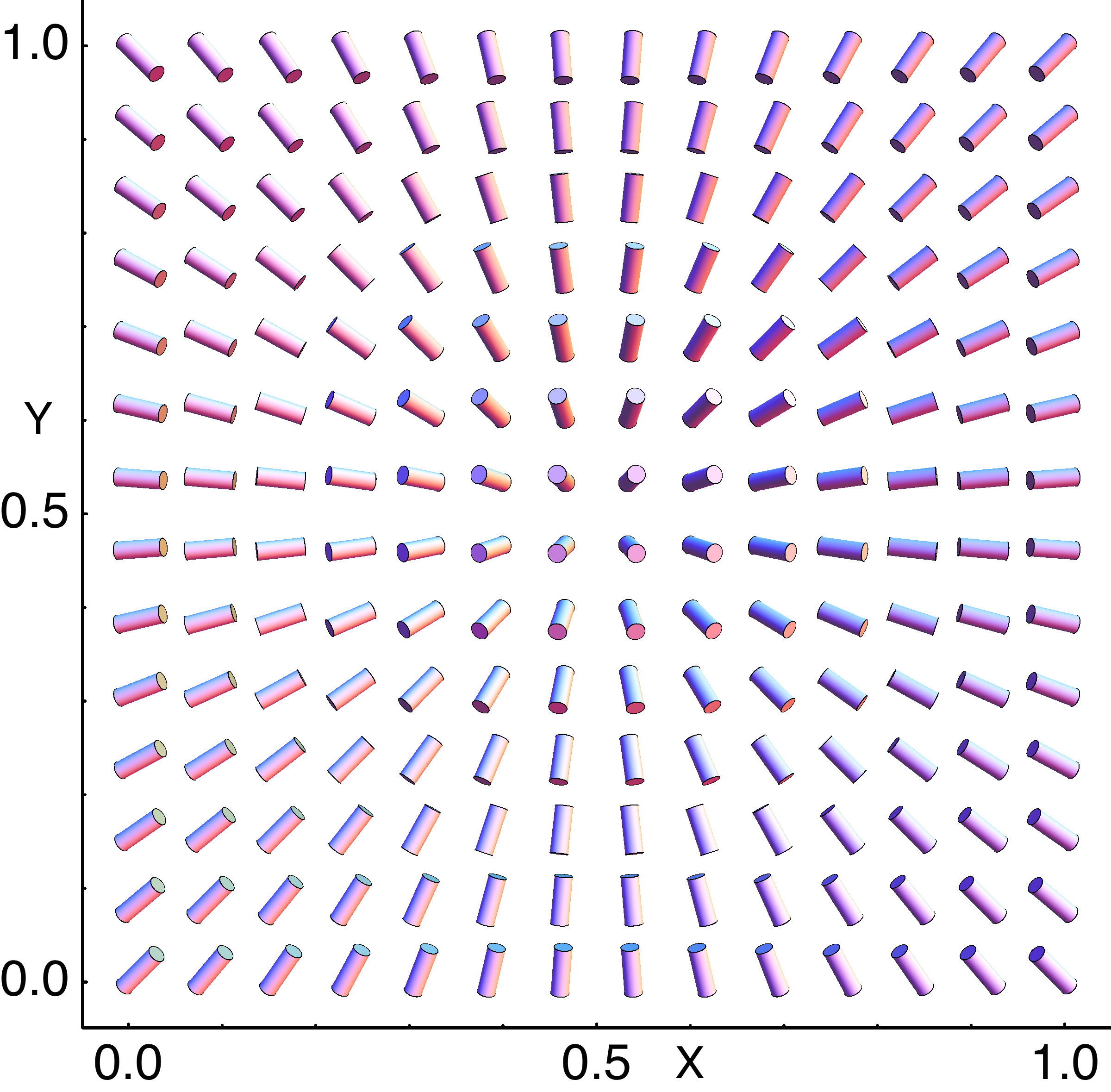} 
      \caption{}
  \label{EscapeDisclinationConfig:right1}
\end{subfigure}
\caption{\small{(\subref{EscapeDisclinationConfig:left1}) Resolved escape solution on $256 \times 256$ mesh (restricted for visualization) with final free energy of $9.971$. (\subref{EscapeDisclinationConfig:center1}) Disclination solution with central wedge defect and final free energy of $24.042$ (free energy is expected to diverge with refinement). (\subref{EscapeDisclinationConfig:right1}) Symmetric escape solution with final free energy of $9.971$.}}
\label{EscapeDisclinationConfig}
\end{figure}

The first solution, located using undeflated solves, is displayed in Figure \ref{EscapeDisclinationConfig}(\subref{EscapeDisclinationConfig:left1}). This director field is continuous and shares some similarities with the solutions found in \cite{Cladis1, Meyer2} for long cylindrical capillaries. The progression of the solves consumes $38.4$ WUs. The solution displayed in Figure \ref{EscapeDisclinationConfig}(\subref{EscapeDisclinationConfig:right1}) is a second, symmetric configuration computed in the deflation solves using $39.5$ WUs. The calculated free energy on each mesh for both solutions is shown in Table \ref{DisclinationFreeEnergies}. Due to the symmetric composition of the device, zenithal tilt in either direction results in an optimal arrangement.

\begin{table}[h!]
\centering
\begin{tabular}[t]{|c|c|c|c|c|c|c|}
\hlinewd{1.3pt}
Grid & $8 \times 8$ & $16 \times 16$ & $32 \times 32$ & $64 \times 64$ & $128 \times 128$ & $256 \times 256$ \\
\hlinewd{1.3pt}
Pos. Escape & $9.972$ & $9.971$ & $9.971$ & $9.971$ & $9.971$ & $9.971$ \\
\hline
Disclination & $13.154$ & $15.331$ & $17.509$ & $19.686$ & $21.864$ & $24.042$ \\
\hline
Neg. Escape & $9.971$ & $9.971$ & $9.971$ & $9.971$ & $9.971$ & $9.971$ \\
\hline
\end{tabular}
\caption{\small{Computed free energies on each mesh for the set of computed solutions.}}
\label{DisclinationFreeEnergies}
\end{table}

In Figure \ref{EscapeDisclinationConfig}(\subref{EscapeDisclinationConfig:center1}), the remaining solution generated through deflation using $121.0$ WUs is displayed. Without deflation, only two of these three solutions would be found across all NI grids. This configuration reveals a disclination where the director field becomes undefined at the center of the domain. The existence of this type of solution lends credence to the escape solution moniker given to configurations like those in Figures \ref{EscapeDisclinationConfig}(\subref{EscapeDisclinationConfig:left1}) and \ref{EscapeDisclinationConfig}(\subref{EscapeDisclinationConfig:right1}), as the director ``escapes'' in the $z$-direction to avoid the defect. Because we assume a slab domain, this is an axial disclination, where the disclination runs parallel to the $z$-axis. The disclination structure does not have finite free energy, as the functional values diverge as they approach the central defect \cite{Stewart1}. Since the solution is approximated with finite elements, this divergent behavior is manifest in a monotonically increasing free energy after each successive refinement; see Table \ref{DisclinationFreeEnergies}. The free energy of $24.042$ computed on the finest grid is expected to continue to rise as the domain is more finely discretized. These types of disclinations can be synthesized and observed under certain conditions \cite{Stephen1}.

\begin{table}[h!]
\centering
\begin{tabular}[t]{|c||c|c|c|c|}
\hline
 & \multicolumn{4}{|c|}{Figure \ref{EscapeDisclinationConfig}} \\
\hlinewd{1.3pt}
Grid & (\subref{EscapeDisclinationConfig:left1}) & (\subref{EscapeDisclinationConfig:center1}) & (\subref{EscapeDisclinationConfig:right1}) & Total Anon. \\
\hlinewd{1.3pt}
$8 \times 8$ & $23$  & $\mathbf{7}$ & $\mathbf{100}$ & $-$ \\
\hline
$16 \times 16$ & $9$ & $12$ & $9$ & $\mathbf{102}$ \\
\hline
$32 \times 32$ & $5$ & $8$ & 5 & $\mathbf{200}$ \\
\hline
$64 \times 64$ & $2$ & $5$ & 2 & $\mathbf{200}$ \\
\hline
$128 \times 128$ & $2$ & $5$ & 2 & $\mathbf{200}$ \\
\hline
$256 \times 256$ & $2$ & $5$ & 2 & $\mathbf{200}$ \\
\hline
\end{tabular}
\caption{\small{Newton iteration counts across grids directly attributable to a solution along with those resulting in divergence or tolerance stoppage in the deflation process for the disclination problem. Bold numbers are associated with the discovery stage using deflation.}}
\label{DisclinationNonlinearIterations}
\end{table}
\vspace{-.2cm}

\begin{table}[h!]
\centering
\begin{tabular}[t]{|c||c|c|c||c|c|c||c|c|c|}
\hline
\ & \multicolumn{3}{|c||}{Fig. \ref{TiltTwistSolutions}: Tilt-Twist} & \multicolumn{3}{|c||}{Fig. \ref{FreederickszTransSolutions}: Freeder.} & \multicolumn{3}{|c|}{Fig. \ref{EscapeDisclinationConfig}: Disclination} \\
\hlinewd{1.3pt}
Grid & (\subref{TiltTwistSolutions:left1}) & (\subref{TiltTwistSolutions:center1}) & (\subref{TiltTwistSolutions:right1}) & (\subref{FreederickszTransSolutions:left1}) & (\subref{FreederickszTransSolutions:center1}) & (\subref{FreederickszTransSolutions:right1}) & (\subref{EscapeDisclinationConfig:left1})  & (\subref{EscapeDisclinationConfig:center1}) & (\subref{EscapeDisclinationConfig:right1}) \\
\hlinewd{1.3pt}
$8 \times 8$& $10.8$ & $\mathbf{11.4}$ & $\mathbf{12.1}$ & $10.7$ & $\mathbf{10.1}$ & $\mathbf{10.1}$ & $12.9$ & $\mathbf{13.3}$ & $\mathbf{13.9}$ \\
\hline
$16 \times 16$ & $8.0$ & $10.5$ & $10.5$ & $10.7$ & $11.0$ & $11.0$ & $13.0$ & $14.9$ & $13.0$ \\
\hline
$32 \times 32$ & $8.0$ & $11.5$ & $11.5$ & $11.0$ & $12.0$ & $12.0$ & $13.0$ & $14.2$ & $13.0$ \\
\hline
$64 \times 64$ & $8.5$ & $12.0$ & $12.0$ & $12.0$ & $12.5$ & $12.5$ & $13.5$ & $16.6$ & $13.5$ \\
\hline
$128 \times 128$ & $9.0$ & $11.5$ & $11.5$ & $12.0$ & $12.5$ & $12.5$ & $14.0$ & $17.8$ & $14.0$ \\
\hline
$256 \times 256$ & $6.0$ & $10.0$ & $10.0$ & $12.0$ & $12.5$ & $12.5$ & $14.0$ & $18.2$ & $14.0$ \\
\hline
\end{tabular}
\caption{\small{Average multigrid iteration counts on each mesh during progression of the NI hierarchy for solutions from the experiments in Section \ref{NumericalExperiments}. Counts in bold represent average iterations for linear solves on deflated systems.}}
\label{NematicMGCounts}
\end{table}

The distribution of nonlinear iterations across the NI
levels for the escape and disclination configurations is shown in
Table \ref{DisclinationNonlinearIterations}. The iteration totals
associated with a particular solution are those Newton steps that
converged to that solution on the given mesh. The counts marked in
bold are iterations performed as part of the discovery stage of the
algorithm using deflation. The total ``anonymous'' iteration counts in
the last column are those deflation steps that resulted in divergence
from the unit-length constraint or reached the Newton iteration limit
without converging. The size of this iteration overhead depends on the
number of initial guesses used and the complexity of the
configurations. While the configurations in this simulation are found
on the first grid, using different deflation parameters can change
this outcome. For instance, using $\alpha = 0.1$ and $p = 1.0$, the
solutions in Figures
\ref{EscapeDisclinationConfig}(\subref{EscapeDisclinationConfig:center1})
and (\subref{EscapeDisclinationConfig:right1}) are not discovered
until the $64 \times 64$ grid. Furthermore, in the next section, a
more complicated example is considered where new solutions are
discovered after significant iteration on finer levels of the hierarchy.

The performance of the linear solver for the three experiments above
is illustrated in Table \ref{NematicMGCounts}. The table displays solver iteration counts averaged over Newton steps on each mesh for the solutions found in each experiment. Bold values delineate average iterations for deflated linear systems. These iteration counts are relatively small and remain steady across mesh refinements. Note that the iteration counts associated with the deflated linear solves are consistent with the performance of the solver on the undeflated systems. These results are especially promising as no special modifications to the solver are necessary for integration with the deflation algorithm.

\section{Cholesteric Liquid Crystals} \label{CholestericNumericalExperiments}

Cholesteric liquid crystals share many properties with nematics but have slightly less symmetry due to chirality. In particular, their inherent helical structure leads to a property known as enantiomorphy where cholesteric molecules are distinguishable from their reflected images. Right-handed helical cholesteric structures are transformed to left-handed helixes upon reflection. This asymmetry leads to a moderate modification of the elastic free-energy functional for these types of liquid crystals and a fourth (nondimensionalized) physical constant, $t_0$. The full free-energy functional is written
\begin{align*}
\mathcal{C}_0(\director) &= K_1 \Ltwoinner{\diverg \director}{\diverg \director}{\Omega} + K_2 \Ltwoinnerndim{\director \cdot \curl \director + t_0}{\director \cdot \curl \director + t_0}{\Omega}{3} \\
&\qquad + K_3 \Ltwoinnerndim{\director \times \curl \director}{\director \times \curl \director}{\Omega}{3} \\
&= \mathcal{F}(\director) + 2K_2 \Ltwoinner{t_0}{\director \cdot \curl \director}{\Omega} + K_2 \Ltwoinner{t_0}{t_0}{\Omega},
\end{align*} 
where $\mathcal{F}(\director)$ is the nematic functional in \eqref{NematicFunctional} without the electric terms. As with the nematic free energies of the previous section, the cholesteric free energies reported here are weighted with the classical factor of $\frac{1}{2}$ for consistency. Note that the last term does not depend on $\director$. Thus, in the minimization process, we need not include that term. Hence, we define the cholesteric free-energy functional to be minimized as
\begin{align*}
\mathcal{C}(\director) = \mathcal{F}(\director) + 2K_2 \Ltwoinner{t_0}{\director \cdot \curl \director}{\Omega}.
\end{align*}

The physical parameter, $t_0$, characterizes the chiral properties of the cholesteric liquid crystal and may be positive or negative depending on the handedness of the cholesteric \cite{Collings1}. Generally, in modeling cholesterics, no saddle-splay term is included regardless of the applied boundary conditions \cite{Stewart1}. While this assumption is sufficient in many practical applications and is used here, cholesteric models relaxing the assumption do exist \cite{Jenkins2}.

\subsection{Minimization}

Since cholesterics are subject to the same pointwise unit-length constraint as nematics, the Lagrangian is formed as
\begin{align*}
\mathcal{L}^C (\director, \lambda) = \mathcal{C}(\director) + \int_{\Omega} \lambda(\director \cdot \director - 1) \diff{V}.
\end{align*}
Computing the derivative of $\mathcal{L}^C$ with respect to $\director$ yields
\begin{align*}
\mathcal{L}_{\director}^C (\director, \lambda)[\vec{v}] = \mathcal{L}_{\director}[\vec{v}] + 2K_2 \left(\Ltwoinner{t_0}{\vec{v} \cdot \curl \director}{\Omega} + \Ltwoinner{t_0}{\director \cdot \curl \vec{v}}{\Omega} \right).
\end{align*}
Since the additional terms of the free energy specific to cholesterics do not depend on $\lambda$, derivatives of this Lagrangian involving $\lambda$ are identical to the nematic case. Thus, in computing the Hessian, the only derivative with additional terms is the second-order derivative with respect to $\director$. This implies that
\begin{align*}
\mathcal{L}^C_{\director \director} = \mathcal{L}_{\director \director} + 2K_2 \left(\Ltwoinner{t_0}{\vec{v} \cdot \curl \ddirector}{\Omega} + \Ltwoinner{t_0}{\ddirector \cdot \curl \vec{v}}{\Omega} \right).
\end{align*}
Modifying the energy-minimization and deflation algorithm discussed
above for nematics by adding in the appropriate terms corresponding to
the cholesteric free energy yields an effective algorithm for
computing multiple equilibrium configurations of cholesteric liquid crystals. 

\subsection{Chiral Configuration} \label{ChiralExample}

The following example is a simple cholesteric configuration that
demonstrates the fundamental departure of cholesteric behavior from
that of nematics. We use the same mixed periodic and Dirichlet
boundary conditions and slab domain assumption as in previous
numerical examples. At the Dirichlet boundary, uniform conditions such
that $\director = (1,0,0)$ are enforced. In the case of nematic liquid
crystals, subject to elastic forces, the minimizing configuration is
full alignment parallel to the director on the boundary. However, the
energetically optimal arrangement for cholesterics is a chiral
configuration along the $y$-axis with twist properties determined by
the value of $t_0$.  Using an ansatz for a chiral solution of the form
$\director = (\cos(\tau y), 0 , -\sin(\tau y))$, the computations in
\cite{Stewart1} can be modified to our coordinate system, giving the
elastic free energy associated with this ansatz as $\frac{1}{2} K_2
(t_0 - \tau)^2 \vert \Omega \vert$, where $\vert \Omega \vert$ is the
domain measure, so long as the chiral ansatz also conforms to the
imposed boundary conditions.  Since the elastic free energy is
positive and semi-definite, clearly the free energy of the ansatz is
minimized when $\tau = t_0$; when $t_0$ is an integer multiple of
$2\pi$, the uniform Dirichlet boundary conditions above will also be satisfied.

For this numerical simulation, the Frank constants are set to $K_1 = 1.0$, $K_2 = 3.0$, and $K_3 = 1.2$ while $t_0 = -2\pi$. This implies that the energy-minimizing solution corresponds to a left-handed helix running parallel to the $y$-axis with a $2\pi$-rotation across the device. However, additional configurations, while not globally minimizing, satisfy the first-order optimality conditions and are experimentally observable.

\begin{figure}[h!]
\centering
\begin{subfigure}{.32 \textwidth}
        \includegraphics[scale=.205, left]{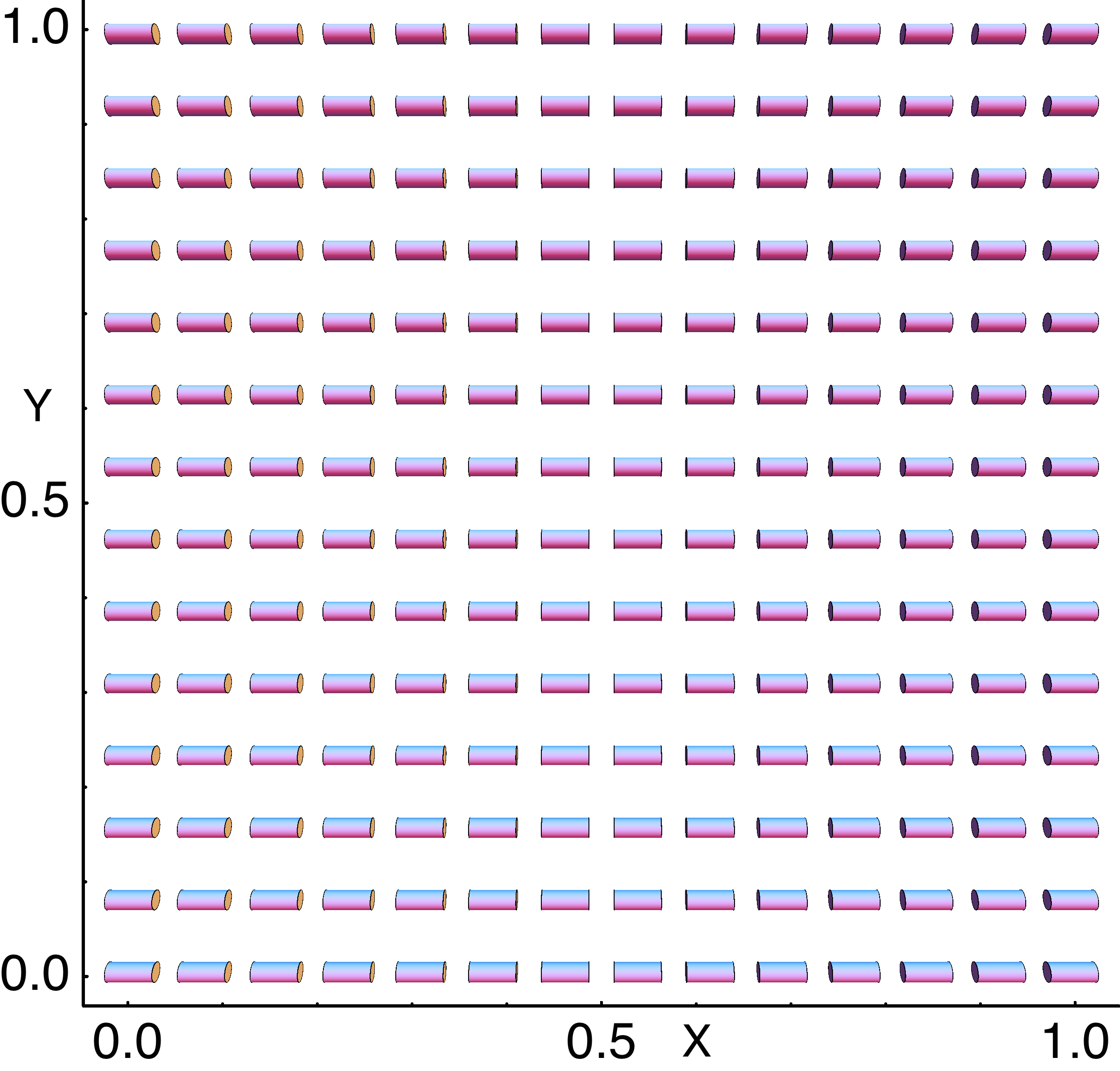}
      \caption{}
  \label{Cholesterics:left1}
\end{subfigure}
\begin{subfigure}{.32 \textwidth} 
  \includegraphics[scale=.205, center]{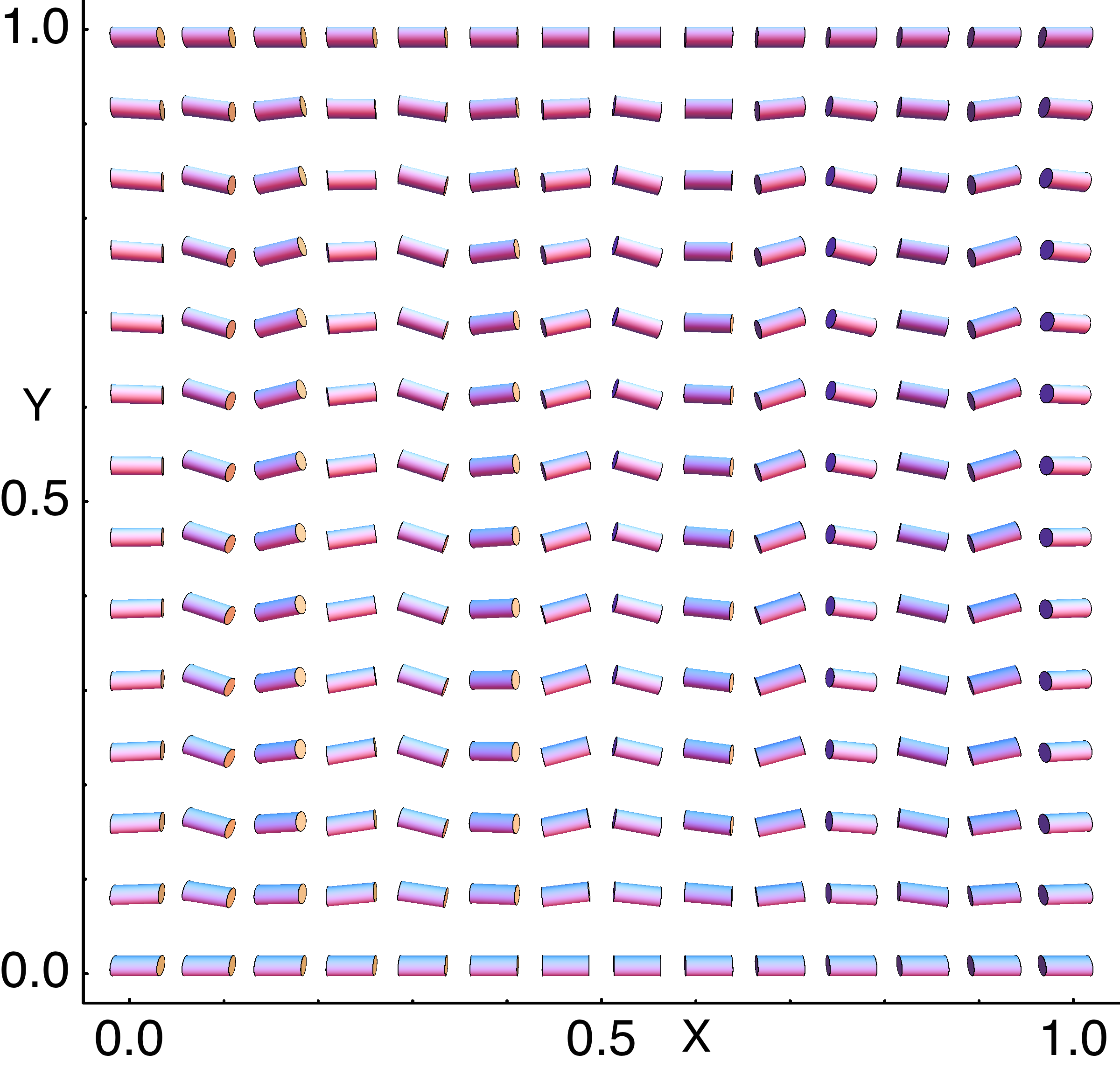}
    \caption{}
  \label{Cholesterics:center1}
\end{subfigure}
\begin{subfigure}{.32 \textwidth}
  \includegraphics[scale=.205, right]{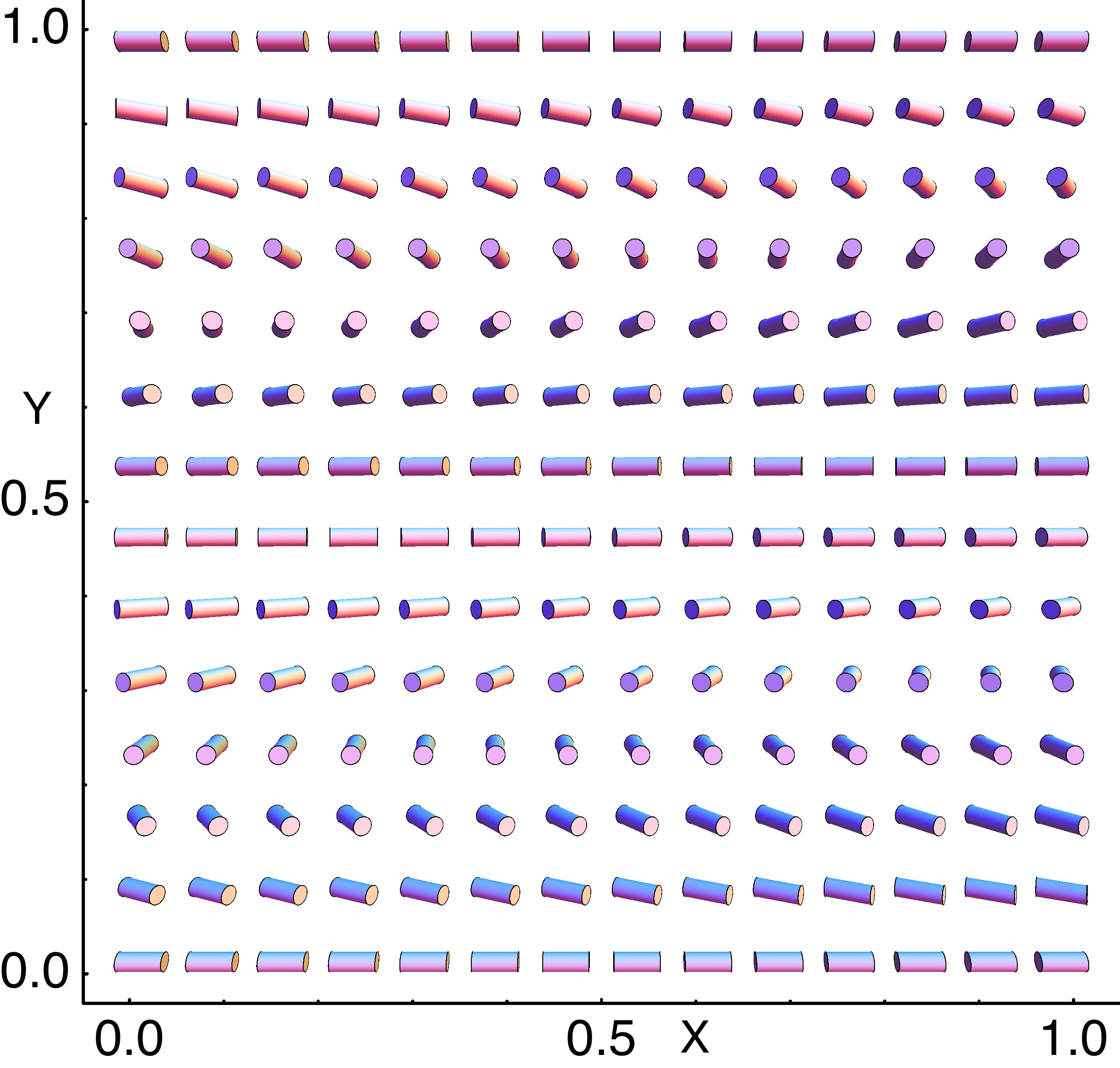} 
      \caption{}
  \label{Cholesterics:right1}
\end{subfigure}
\begin{subfigure}{.32 \textwidth}
  \includegraphics[scale=.205, left]{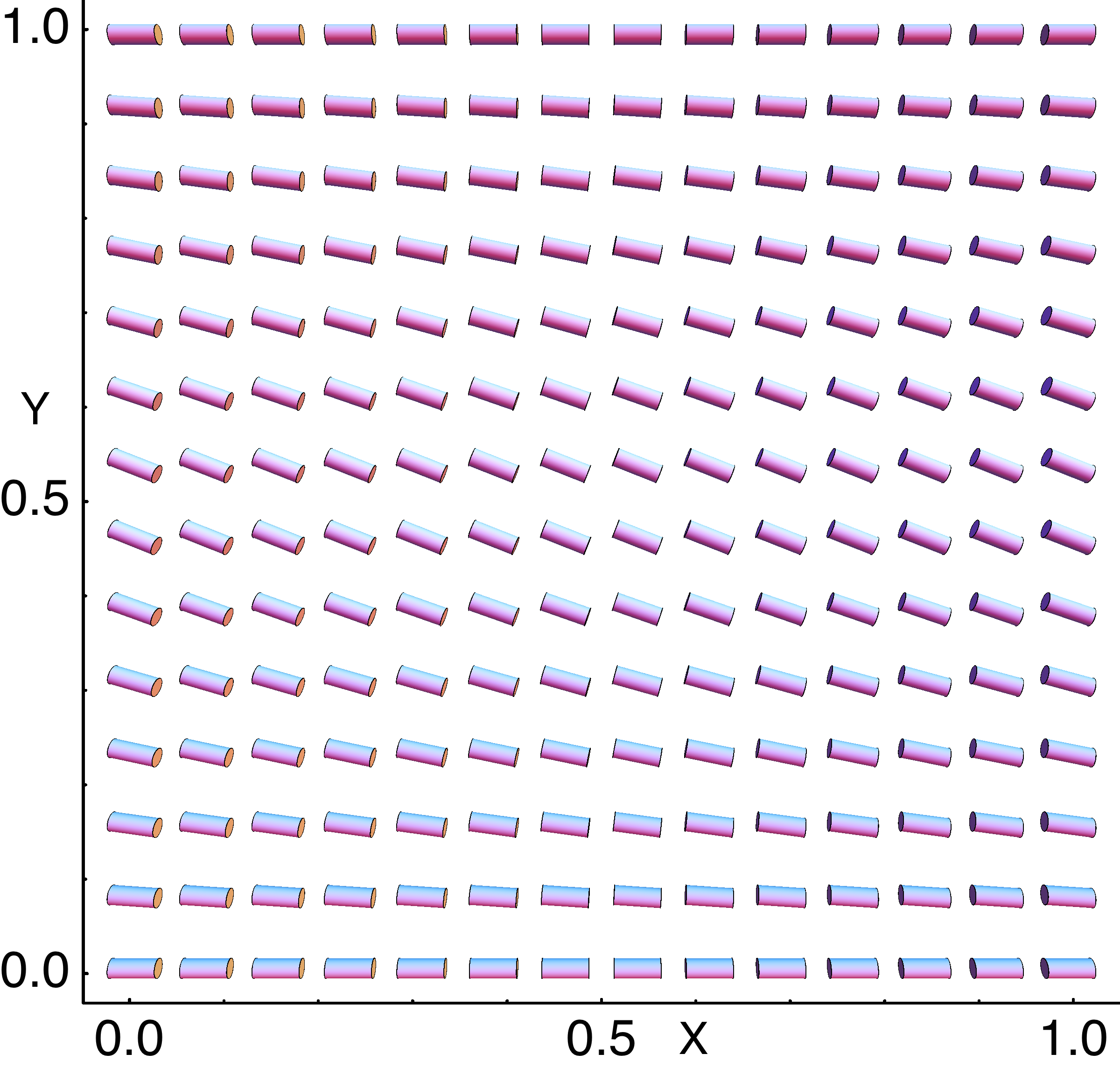}
      \caption{}
  \label{Cholesterics:left2}
\end{subfigure}
\begin{subfigure}{.32 \textwidth}
  \includegraphics[scale=.205, center]{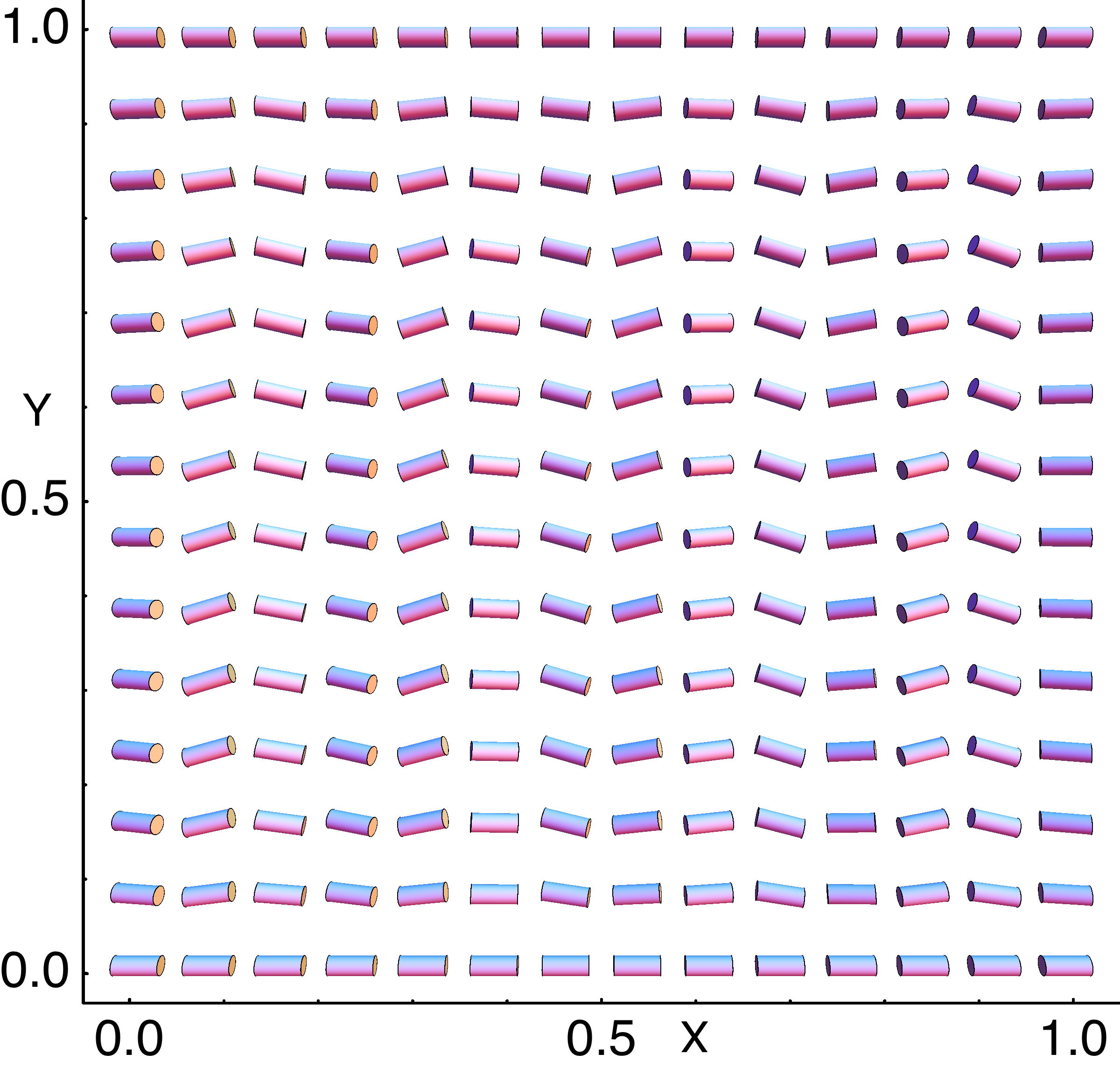} 
      \caption{}
  \label{Cholesterics:center2}
\end{subfigure}
\begin{subfigure}{.32 \textwidth}
  \includegraphics[scale=.205, right]{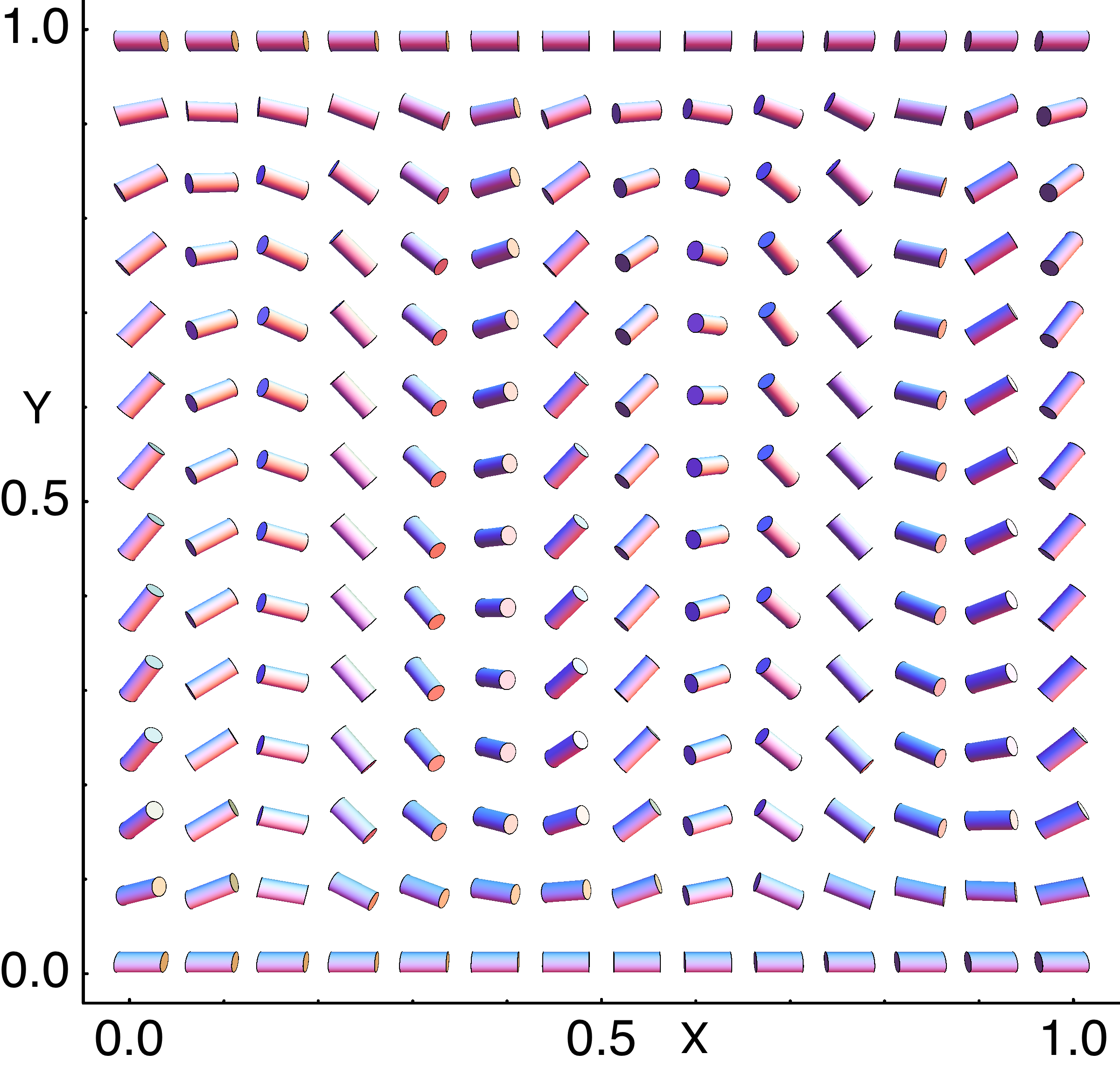}
      \caption{}
  \label{Cholesterics:right2}
\end{subfigure}
\caption{\small{Family of distinct solutions for the cholesteric equilibrium problem found through deflation. Each solution is computed computed on a $256 \times 256$ mesh and restricted for visualization. The energy minimizing solution is displayed in (\subref{Cholesterics:right1}).}}
\label{CholestericSolutionConfig}
\end{figure}

The deflation algorithm is applied with damping values of $\omega_1 = 0.2$, $\Delta_1 = 0.2$, $\omega_2 = 0.2$, and $\Delta_2 = 0.0$. Using the set of three initial guesses outlined in Appendix \ref{InitialGuesses}, the algorithm reveals a rich set of solutions satisfying the optimality conditions. A total of six distinct solutions, shown in Figure \ref{CholestericSolutionConfig}, are found, whereas, without deflation, only three solutions would be identified across all grids. The corresponding, computed free energies for these solutions is shown in Table \ref{CholestericMGCounts} along with average iteration counts for the multigrid-preconditioned GMRES linear solver. In general, the linear solver iteration counts are higher for these cholesteric systems compared with those of the previous section. Correspondingly, the WUs, shown in the same table, are larger when compared with previous experiments. This increase in iterations is most likely due to a combination of the additional term in the cholesteric functional and higher overall free energies in the solutions. However, the iterations counts are relatively consistent across grid refinements, and the average solver iterations for the deflated systems, shown in bold, correspond well with the iteration counts for the undeflated solves.

The solution set includes degenerate planar solutions displayed in Figures \ref{CholestericSolutionConfig}(\subref{Cholesterics:left1}) and (\subref{Cholesterics:left2}). By virtue of the chiral proclivity of cholesterics, these configurations are not globally minimizing. Also included in the computed arrangements are structures with left-hand twists following the $x$-axis. These twisting behaviors decrease each system's free energy below that of the planar solutions. For instance, the free energy of the configuration in Figure \ref{CholestericSolutionConfig}(\subref{Cholesterics:right2}) is $31.821$, well below the free energy of $59.218$ for the solution in Figure \ref{CholestericSolutionConfig}(\subref{Cholesterics:left1}). These transversal wave-like structures can be observed in certain cholesteric samples. Note that this solution is not located until the fifth mesh in the NI framework. Furthermore, if the deflation parameters are adjusted, an additional configuration of this type, not shown here, is located during the deflation process.

As in Table \ref{DisclinationNonlinearIterations}, the nonlinear iterations accrued on each mesh during the NI process are presented in Table \ref{CholestericNonlinearIterations}. Iteration counts in bold delineate those performed using deflation. Note that the deflation iteration counts are generally higher with the increased Newton damping in the deflation stage. Further, distinct solutions are discovered through deflation on multiple grids. This is especially apparent in the discovery of the structure in Figure \ref{CholestericSolutionConfig}(\subref{Cholesterics:right2}), which is only found on the second finest mesh. Thus, eliminating or reducing deflation iterations on finer grids is difficult without risking the loss of additional solutions.

\begin{table}[h!]
\centering
\begin{tabular}[t]{|c||c|c|c|c|c|c|c|c|c|c|c|}
\hline
 & \multicolumn{6}{|c|}{Figure \ref{CholestericSolutionConfig}} \\
\hlinewd{1.3pt}
Grid &(\subref{Cholesterics:left1}) & (\subref{Cholesterics:center1}) & (\subref{Cholesterics:right1}) & (\subref{Cholesterics:left2}) & (\subref{Cholesterics:center2}) & (\subref{Cholesterics:right2}) \\
\hlinewd{1.3pt}
$8 \times 8$ & $46.2$ & $\mathbf{52.9}$ & $\mathbf{11.3}$ & $-$ & $-$ & $-$ \\
\hline
$16 \times 16$ & $66.0$ & $54.3$ &  $9.0$ & $\mathbf{67.7}$ & $\mathbf{66.8}$ & $-$ \\
\hline
$32 \times 32$ & $65.0$ & $33.2$ & $8.0$ & $53.9$ & $34.1$ & $-$ \\
\hline
$64 \times 64$ & $61.0$ & $28.4$ & $8.0$ & $35.8$ & $26.9$ & $-$ \\
\hline
$128 \times 128$ & $62.0$ & $33.0$ & $9.0$ & $52.5$ & $32.5$ & $\mathbf{29.0}$ \\
\hline
$256 \times 256$ & $78.0$ & $30.5$ & $9.5$ & $46.0$ & $30.0$ & $18.5$ \\
\hlinewd{1.3pt}
Work Units & $100.7$ & $103.7$ & $28.5$ & $156.9$ & $108.7$ & $493.0$ \\
\hline
Free Energy & $59.218$ & $56.553$ & $2.984$e-$08$ & $59.378$ & $56.553$ & $31.821$ \\
\hline
\end{tabular}
\caption{\small{Average multigrid iteration counts on each mesh during progression of the NI hierarchy for selected solutions from the cholesteric experiment above. Counts in bold represent average iterations for linear solves on deflated systems. The final rows display the WUs and free energy associated with each computed equilibrium configuration.}}
\label{CholestericMGCounts}
\end{table}

\begin{table}[h!]
\centering
\begin{tabular}[t]{|c||c|c|c|c|c|c|c|}
\hline
 & \multicolumn{7}{|c|}{Figure \ref{CholestericSolutionConfig}} \\
\hlinewd{1.3pt}
Grid & (\subref{Cholesterics:left1}) & (\subref{Cholesterics:center1}) & (\subref{Cholesterics:right1}) & (\subref{Cholesterics:left2}) & (\subref{Cholesterics:center2}) & (\subref{Cholesterics:right2}) & Total Anon. \\
\hlinewd{1.3pt}
$8 \times 8$ & $46$ & $\mathbf{56}$ & $\mathbf{50}$ & $-$ & $-$ & $-$ & $\mathbf{100}$ \\
\hline
$16 \times 16$ & $1$ & $22 $ & $19$ & $\mathbf{87}$ & $\mathbf{55}$ & $-$ & $\mathbf{100}$ \\
\hline
$32 \times 32$ & $1$ & $12$ & $10$ & $8$ & $12$ & $-$ & $\mathbf{228}$ \\
\hline
$64 \times 64$ & $1$ & $7$ & $5$ & $4$ & $7$ & $-$ & $\mathbf{233}$ \\
\hline
$128 \times 128$ & $1$ & $2$ & $2$ & $2$ & $2$ & $\mathbf{63}$ & $\mathbf{200}$  \\
\hline
$256 \times 256$ & $1$ & $2$ & $2$ & $2$ & $2$ & $2$ & $\mathbf{253}$ \\
\hline
\end{tabular}
\caption{\small{Newton iteration counts across grids directly attributable to a solution along with those resulting in divergence or tolerance stoppage in the deflation process for the cholesteric problem. Bold numbers are associated with the discovery stage using deflation.}}
\label{CholestericNonlinearIterations}
\end{table}

The arrangement displayed in Figure \ref{CholestericSolutionConfig}(\subref{Cholesterics:right1}), is the energetically optimal configuration. As suggested in the analysis above, the director profile contains a left-handed helical structure rotating $2\pi$ radians about an axis parallel to the $y$-axis. In agreement with the derived analytical free energy, the computed free energy for this solution is $2.984$e-$08$. It should be noted that in order to obtain this solution, an initial guess incorporating a twisting profile is used. Moving in the configuration space from a profile with little or no twist to one that incorporates a full $2 \pi$-rotation is far from monotonic in terms of energy optimization. For example, introducing a moderate twist into the planar solutions of Figures \ref{CholestericSolutionConfig}(\subref{Cholesterics:left1}) or (\subref{Cholesterics:left2}) increases their free energy until the twist approaches a $2 \pi$-rotation. The minima valleys are well-separated and the strength of the poles introduced through deflation are often not enough to overcome the barrier dividing the valleys. Investigation into the application of generalized tunneling methods, which have been used to address some aspects of this challenge in the context of function minimization \cite{Levy1}, will be the subject of future work.

\section{Conclusion and Future Work} \label{Conclusions}

We have discussed a deflation technique for the computation of distinct solutions in the context of a free-energy variational approach for the simulation of nematic and cholesteric liquid crystal equilibrium configurations under the Frank-Oseen model. It was shown that highly accurate and efficient multigrid methods previously designed for the original undeflated discrete systems are applicable in solving the associated deflated linear systems. To further increase the efficacy and efficiency of the approach, a strategy for interweaving the deflation technique with an NI framework was presented. This produces a dynamic algorithm and enables the discovery of additional equilibrium configurations. 

Four illustrative numerical simulations were conducted with the proposed algorithm. These results demonstrate the power of deflation to systematically resolve  complicated bifurcation phenomena and disclination behavior in nematics, as well as chiral configurations in cholesterics. In each application, deflation locates multiple configurations satisfying the first-order optimality conditions, including both local and global extrema. Each simulation successfully employed a coupled multigrid method with Braess-Sarazin-type relaxation and NI to improve overall efficiency. Future work will consider construction of a generalized tunneling approach, based on the work in \cite{Levy1}, applied to the Newton iterations to further increase the power of the deflation method. In addition, we aim to investigate the method's performance in analyzing new physical phenomena and behaviors in shaped domains. Finally, strategies for dynamic and adaptive construction of initial guesses for the deflation method at each level of the NI process will be studied.

\section*{Acknowledgements}

The authors would like to thank Professor Timothy Atherton for his
useful suggestions and guidance. We would also like to thank
Dr. Thomas Benson for allowing us to adapt his code.

\bibliographystyle{siam}

\bibliography{MathematicalCitations}

\appendix

\section{Initial Guesses} \label{InitialGuesses}

In this appendix we report the initial guesses used for each example,
to aid in reproducing the results. Each guess listed here gives the values used on the interior of the domain for all NI levels; the Dirichlet boundary conditions are enforced along the relevant boundaries. In all of the simulations performed, $\lambda$ is initially set to $0$. 

In Sections \ref{TiltTwistExperiments} and
\ref{FreederickszExperiments}, the initial guesses used were
$\director = \left (\cos \left(\frac{\pi}{40} \right), \sin
  \left(\frac{\pi}{40} \right), 0 \right)$ and $\director = \left
  (\cos \left(\frac{\pi}{40}\right), -\sin \left(\frac{\pi}{40}
  \right), 0 \right)$. In addition, the simulations of Section
\ref{FreederickszExperiments} use $\phi = V \cdot y$ to initialize the
electric potential for both guesses, where $V$ is the potential at the
top substrate. 

For Section \ref{EscapeDisclination}, let $\xi_1 = \left \vert \tan^{-1} \left (\frac{0.5 - y}{0.5 - x} \right) \right \vert$, $\zeta_1 = \frac{9 \pi}{20}$, and define the functions
\begin{align*}
n_1 &= \begin{cases}
\sin(\zeta_1) \cos(\xi_1) & \text{if $x \leq 0.5$} \\
-\sin(\zeta_1) \cos(\xi_1) &  \text{if $x > 0.5$},
\end{cases} &
n_2 &= \begin{cases}
\sin(\zeta_1) \sin(\xi_1) & \text{if $y \leq 0.5$} \\
-\sin(\zeta_1) \sin(\xi_1) & \text{if $y > 0.5$},
\end{cases} \\
n_3 &= \cos(\zeta_1).
\end{align*}
Then the two initial values for the director in the section are given by
\begin{align*}
\director &= \begin{cases}
\big( 0, 0, 1 \big) & \text{if $x, y = 0.5$} \\
\big( n_1, n_2, n_3 \big) & \text{otherwise}, \\
\end{cases} & 
\director &= \begin{cases}
\big( 0, 0, 1 \big) & \text{if $x, y = 0.5$} \\
\big( n_1, n_2, -n_3 \big) & \text{otherwise}. \\
\end{cases} 
\end{align*}

Finally, for Section \ref{ChiralExample}, let $\xi_2 = \frac{7 \pi}{16}$ and $\zeta_2 = \frac{\pi}{4}$. The initial values for $\director$ are shown in Table \ref{InitialGuessFormulas}.

\begin{table}[h!]
\centering 
\small{
\begin{tabular}[t]{|c|c|c|}
\hline
Guess $1$ & Guess $2$ & Guess $3$ \\
\hlinewd{1.3pt} 
$\begin{aligned}
n_1 &= \cos \left(\pi/12 \right) \\
n_2 &= \sin \left(\pi/12 \right) \\
n_3 &= 0 \\
\end{aligned}$ & 
$\begin{aligned}
n_1 &= \sin \left(\xi_2 \right) \cos \big( \zeta_2 \cos(4 \pi x) \big)\\
n_2 &= \sin \left(\xi_2 \right) \sin \big( \zeta_2 \cos(4 \pi x) \big) \\
n_3 &= \cos(\xi_2) \\
\end{aligned}$ &
$\begin{aligned}
n_1 &= \cos \left(2 \pi y \right) \cos\left( \pi/8 \right) \\
n_2 &= \cos \left(2 \pi y \right) \sin \left(\pi/8 \right) \\
n_3 &= \sin (2 \pi y) \\
\end{aligned}$ \\
\hline
\end{tabular}}
\caption{\small{Formulas for the initial guesses used in Section \ref{ChiralExample}.}}
\label{InitialGuessFormulas}
\end{table}

\end{document}